\documentclass[titlepage,11pt]{article}
\usepackage{latexsym}
\usepackage{amsfonts}
\usepackage{amsmath}
\usepackage{tikz}
\usepackage{float}
\newcommand\blackslug{\hbox{\hskip 1pt \vrule width 4pt height 8pt depth 1.5pt
        \hskip 1pt}}
\newcommand\bbox{\hfill \quad \blackslug \bigbreak}
\def\DD{\hbox{-}}

\def\LL{,\ldots,}
\def\cupcup{\cup\cdots\cup}
\newcommand{\vare}{\varepsilon}

\title{Bipartite graphs with no $K_6$ minor}
\author{
Maria Chudnovsky\thanks{Supported by NSF DMS-EPSRC grant DMS-2120644.}\\
Princeton University, Princeton, NJ 08544
\\
\\
Alex Scott\thanks{Research supported by EPSRC grant EP/V007327/1.}\\
Mathematical Institute, University of Oxford, Oxford OX2 6GG, UK
\\
\\
Paul Seymour\thanks{Supported by AFOSR grant
A9550-19-1-0187, and by NSF grant  DMS-1800053.}\\
Princeton University, Princeton, NJ 08544
\\
\\
Sophie Spirkl\thanks{We acknowledge the support of the Natural Sciences and Engineering Research
Council of Canada (NSERC), [funding reference number RGPIN-2020-03912].Cette
recherche a \'et\'e financ\'ee par le Conseil de recherches en sciences naturelles et
en g\'enie du Canada (CRSNG), [num\'ero de r\'ef\'erence RGPIN-2020-03912].}\\
University of Waterloo, Waterloo, Ontario N2L3G1, Canada}

\date{February 19 2022; revised \today}

\newtheorem{thm}{}[section]

\newcommand{\Proof}{\noindent{\bf Proof.}\ \ }

\begin{document}
\maketitle
\begin{abstract}
A theorem of Mader shows that every graph with average degree at least eight has a $K_6$ minor, and this is false if we replace 
eight by any smaller constant. Replacing average degree by minimum degree seems to make little difference:  we do not know 
whether all graphs with minimum degree at least seven have $K_6$ minors, but minimum degree six is certainly not enough.
For every $\vare>0$ there are 
arbitrarily large graphs with average degree at least $8-\vare$ and minimum degree at least six, with no $K_6$ minor.

But what if we restrict ourselves to bipartite graphs? The first statement remains true:
for every $\vare>0$ there are 
arbitrarily large bipartite graphs with average degree at least $8-\vare$ 
and no $K_6$ minor. But surprisingly, going to minimum degree now makes a significant difference.
We will show 
that every bipartite graph with minimum degree at least six
has a $K_6$ minor. Indeed, it is enough that every vertex in the larger part of the bipartition has degree at least six.

\end{abstract}

\section{Introduction}
The graphs with no $K_5$ minor are well understood. A theorem of Wagner~\cite{wagner} gives an explicit construction for all such 
graphs: they can all be built by piecing together planar graphs and copies of one eight-vertex graph by a sum operation that we 
do not describe here. Consequently, every graph with $n\ge 3$ vertices and more than $3n-6$ edges has a $K_5$ minor. 
(All graphs in this paper are finite and have no loops or parallel edges.) This is tight: there are graphs with $n$ vertices 
and with exactly $3n-6$ edges that have no $K_5$ minor. Indeed, one can make such graphs that are almost $6$-regular: for infinitely 
many values of $n$ there is an $n$-vertex planar graph (which therefore has no $K_5$ minor) with all vertices of degree six
except for twelve of degree five. In summary: 
\begin{itemize}
\item all graphs with average degree at least six contain $K_5$ minors, and this is false if we 
replace six by any smaller real number; 
\item all graphs with minimum degree at least six have $K_5$ minors, and this is false if 
we replace six by any smaller integer; 
\item this is all still true even if we insist that maximum degree is at most six.
\end{itemize}

What if we look just at bipartite graphs? One can make $n$-vertex bipartite graphs with no $K_5$ minor that have $3n-9$ edges
(the complete bipartite graph $K_{3,n-3}$ -- in fact this is the only such graph, which can easily be shown by induction using
\ref{mindeg4}). So the situation for average degree is virtually unchanged: average degree six is 
enough to guarantee a $K_6$ minor, and no smaller constant works. 
But $K_{3,n-3}$
has vertices with
degree much larger than the average, and also vertices with degree much smaller than average (if three is much smaller than six).
So what happens if we insist that maximum degree is close to the average degree, or minimum degree is large?

It turns out that:
\begin{thm}\label{mindeg4}
Every non-null bipartite graph with minimum degree at least four has a $K_5$ minor. 
\end{thm}
This can be derived from
Wagner's construction~\cite{wagner}, although the proof is rather long and we omit it. The result is already known: it was stated 
(in a stronger form, replacing ``bipartite'' by ``girth at least four'') in a lecture by
J\'{a}nos Bar\'{a}t~\cite{barat}, as joint work with David Wood, and also (without proof) in an early version of the 
paper~\cite{barat1} (unfortunately it was removed in a later version of the paper).
When excluding $K_5$, imposing a bound on 
maximum degree is perhaps not so interesting: there are $n$-vertex bipartite graphs with average degree at least four and 
maximum degree at most five, with no $K_5$ minor. (For example, take five disjoint copies of $K_{3,5}$, and for $1\le i\le 5$
let $v_i$ be a vertex with degree three from the $i$th copy. Now add two more vertices both adjacent to each of $v_1\LL v_5$. 
To make bigger examples, take disjoint unions.)
Perhaps average degree at least five and maximum degree at most six will guarantee a $K_5$ minor in a bipartite graph, but 
we have not worked this out.


In this paper, we ask what happens for $K_6$ minors. A theorem of Mader~\cite{mader} says:
\begin{thm}\label{mader}
For $n\ge 4$, every $n$-vertex graph with more 
than $4n-10$ edges has a $K_6$ minor.
\end{thm}
There are graphs with $n$ vertices and $4n-10$ edges with minimum degree at least six that have 
no $K_6$ minor: for instance, take a planar graph on $n-1$ vertices with $3(n-1)-6$ edges and minimum degree five, and add a new 
vertex adjacent to everything. So we need average degree at least eight to guarantee a $K_6$-minor;
no smaller constant works. Again, we might ask what happens if we insist that maximum degree is close to average degree, or
minimum degree is large. We have found $K_6$-minor-free graphs with minimum degree six and maximum degree at most nine; and $K_6$-minor-free graphs 
with minimum degree five, maximum degree seven, and average degree arbitrarily close to
$98/15$ (we omit the details). But as far as we know, both the following are open:
\begin{thm}\label{6-regular}
{\bf Conjecture: } Every non-null 6-regular graph has a $K_6$ minor.
\end{thm}
(Indeed, as far as we know, every non-null graph with minimum degree at least six and maximum degree at most eight has a $K_6$ minor.)
\begin{thm}\label{gen6}
{\bf Conjecture: } Every non-null graph with minimum degree at least seven has a $K_6$ minor.
\end{thm}

There is a well-known conjecture of J{\o}rgensen~\cite{jorgensen} that is related:
\begin{thm}\label{jorgconj}
{\bf Conjecture: } Every 6-connected graph with no $K_6$ minor can be made planar by deleting some vertex, and therefore has a vertex of degree at most six.
\end{thm}

But in this paper we will restrict ourselves to bipartite graphs.
Still no constant smaller than eight works as a bound on
average degree to guarantee a $K_6$ minor, since the complete bipartite graph $K_{4,n-4}$ has $4n-16$ edges and has no $K_6$ minor.
But what about minimum degree? We will show:
\begin{thm}\label{symthm}
Every non-null bipartite graph with minimum degree at least six has a $K_6$ minor.
\end{thm}

We do not know whether ``six'' can be replaced by ``five'' in \ref{symthm}. Minimum degree is more difficult than average degree 
to work with inductively, and fortunately there is a strengthening of \ref{symthm} that is more amenable to induction:

\begin{thm}\label{mainthm}
Let $G$ admit a bipartition $(A,B)$ with $|A|\ge |B|>0$, such that every vertex in $A$ has degree at least six. Then $G$ has a $K_6$
minor.
\end{thm}
We remark that \ref{mainthm} becomes false if we replace ``six'' by ``five''; we will show this in the next section.

This was also motivated by one of the steps in the proof of~\cite{RST} that every graph with no $K_6$ minor is five-colourable.
Let $G$ be a minor-minimal graph with no $K_6$ minor that is not five-colourable, if such a graph exists; then in~\cite{RST}, 
section 12 was devoted to showing that
$G$ has a matching with at least $(|G|-1)/2$ edges. If not, then by Tutte's theorem, there is a set $X\subseteq V(G)$
such that $G\setminus X$ has more than $|X|$ odd components, and it was known that $G$ is six-connected, and so each of these
components has an edge to at least six vertices in $X$. By contracting these components to single vertices we obtain a
bipartite graph satisfying the hypotheses of \ref{mainthm}, which would be a contradiction, since $G$ has no $K_6$ minor.
In~\cite{RST}, Mader's theorem~\cite{mader} was used in place of \ref{mainthm}, with additional analysis of the components that had only six 
or seven neighbours in $X$. But it should be added that \ref{mainthm} is not going to shorten the proof of the main theorem of~\cite{RST};
the proof of \ref{mainthm} is considerable longer than section 12 of~\cite{RST}.
It will take up almost all the paper, but we begin with proving the statements for $K_5$ mentioned above.

\section{Some definitions, and the results for $K_5$}
Let us be more precise. If $X\subseteq V(G)$, $G\setminus X$ is the graph obtained from $G$ by deleting $X$, and 
$G[X]=G\setminus (V(G)\setminus X)$ denotes the subgraph of $G$ induced on $X$. A graph $H$ is a {\em minor} of $G$ if $H$ can be obtained by edge-contraction from a subgraph of $G$.
(We repeat that graphs in this paper have no loops or parallel edges, so any loops or parallel edges produced by edge-contraction
should be deleted.) We will only be concerned with complete graph minors. Let us say a {\em cluster} in $G$ is a set of 
disjoint subsets $X_1\LL X_k$ of $V(G)$, such that $G[X_i]$ is connected for $1\le i\le k$, and for $1\le i<j\le k$ there is an edge
of $G$ between $X_i, X_j$; and a {\em $t$-cluster} means a cluster of cardinality $t$.
Thus $G$ contains the complete graph $K_t$ as a minor if and only if $G$ contains a $t$-cluster.

A word on taking minors of bipartite graphs: we start with a graph with a bipartition $(A,B)$, choose a subset $X\subseteq V(G)$ that induces a connected subgraph, and contract $X$ to a single vertex. As we said, if this produces parallel edges we delete them, since we only work 
with simple graphs in this paper. But there is another issue: the graph we obtain by contraction might not be bipartite, and we want 
to produce a bipartite graph at the end, so we in general we must delete some of the edges incident with the new vertex. We could 
explicitly list the edges that we need to delete, but since we will apply this operation many times, let us set up a more convenient 
method. Let us say we {\em contract $X$ into $A$} if we first contract $X$ to a single vertex, $x$ say, and then delete all edges between $x$
and $A$. Thus the graph we produce has a bipartition $((A\setminus X)\cup \{x\},B\setminus X)$. ``Contracting into $B$'' is defined 
similarly.

Let us see first: 
\begin{thm}\label{smallcases}
For $t=1,2,3,4$, if $G$ admits a bipartition $(A,B)$ with $|A|\ge |B|>0$ such that every vertex in $A$ has degree at least $t-1$
then $G$ has a $K_t$ minor.
\end{thm}
\Proof We may assume that every vertex in $A$ has degree exactly $t-1$, by deleting edges, and we may assume that
$|A|=|B|$, by deleting $|A|-|B|$ vertices from $A$.
For $t\le 2$ the result is clear. For $t=3$, the graph has $2|A| = |G|$ edges and so has a cycle, and hence 
a $K_3$ minor. 

Next let $t=4$; we proceed by induction on $|A|$. We may assume that $G$ has a vertex
of degree at most two, $b$ say (necessarily $b\in B$), because otherwise it has a $K_4$ minor. 
If $b$ has degree zero we may delete it, and if it has degree one we may 
delete it and its neighbour, and in either case the result follows from the inductive hypothesis. So we assume that
$b$ has two neighbours 
$a_1,a_2$. If there are at least four vertices in $B\setminus \{b\}$ with a neighbour in $\{a_1,a_2\}$, we may contract 
$\{a_1,b,a_2\}$ 
into $A$
and apply the inductive hypothesis; so we assume that $a_1,a_2$ have exactly the same neighbours $b,b_1,b_2,b_3$.
If some vertex different from $a_1,a_2$ is adjacent to both $b_1,b_2$ then $G$ has a $K_4$ minor: and otherwise
we may contract $\{b_1,b_2,a_1,a_2,b\}$ into $B$ and apply the inductive hypothesis.
This proves \ref{smallcases}.~\bbox

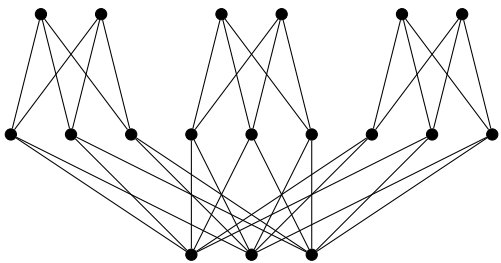
\begin{figure}[H]
\centering

\begin{tikzpicture}[scale=0.8,auto=left]
\tikzstyle{every node}=[inner sep=1.5pt, fill=black,circle,draw]
\node (a1) at (-4,2) {};
\node (a2) at (-3,2) {};
\node (a3) at (-2,2) {};

\node (a4) at (-1,2) {};
\node (a5) at (0,2) {};
\node (a6) at (1,2) {};

\node (a7) at (2,2) {};
\node (a8) at (3,2) {};
\node (a9) at (4,2) {};
\node (b1) at (-1,0) {};
\node (b2) at (0,0) {};
\node (b3) at (1,0) {};

\node (b4) at (-3.5,4) {};
\node (b5) at (-2.5,4) {};
\node (b6) at (-.5,4) {};
\node (b7) at (.5,4) {};
\node (b8) at (2.5,4) {};
\node (b9) at (3.5,4) {};

\foreach \from/\to in {b1/a1,b1/a2,b1/a4, b1/a5,b1/a7,b1/a8,b2/a1,b2/a3,b2/a4,b2/a6,b2/a7,b2/a9,b3/a2,b3/a3,b3/a5,b3/a6,b3/a8,b3/a9,a1/b4,a1/b5,a2/b4,a2/b5,a3/b4,a3/b5,a4/b6,a4/b7,a5/b6,a5/b7,a6/b6,a6/b7,a7/b8,a7/b9,a8/b8,a8/b9,a9/b8,a9/b9}
\draw [-] (\from) -- (\to);

\end{tikzpicture}
\caption{Counterexample to \ref{smallcases} with $t=5$.} \label{fig:smallcases}
\end{figure}

Since all bipartite graphs with minimum degree at least four have $K_5$ minors, 
one might hope that \ref{smallcases} would hold with $t=5$, but that is false. Here is a counterexample 
(see figure \ref{fig:smallcases}). Let $H$ be the graph 
obtained from $K_{3,5}$ by deleting three edges that form a matching. 
Now take $k$ copies of $H$, say $H_1\LL H_k$, and for $1\le i\le k$ let $a_i,b_i,c_i$ be the three vertices of $H_i$ that have degree two.
Let $G$ be obtained from the disjoint union of $H_1\LL H_k$ by making the identifications $a_1=\cdots=a_k$,
$b_1=\cdots=b_k$ and $c_1=\cdots= c_k$. Then $G$ admits a bipartition $(A,B)$ with $|A|=3k$ and $|B|=2k+3$, and every vertex 
in $A$ has degree four, and $G$ has no $K_5$
minor. Thus taking $k\ge 3$ we obtain a counterexample to \ref{smallcases} with $t=5$. By taking $k=4$ instead, and then adding
a new vertex adjacent to every vertex in $A$, we obtain a graph that shows that we cannot replace ``six'' by ``five'' in \ref{mainthm}.

We have:
\begin{thm}\label{smallcases5}
If $G$ admits a bipartition $(A,B)$ with $|A|\ge |B|>0$ such that every vertex in $A$ has degree at least five
then $G$ has a $K_5$ minor.
\end{thm}
This is turn is a consequence of \ref{mainthm} as we show now.

\bigskip

\noindent{\bf Proof of \ref{smallcases5}, assuming \ref{mainthm}. \ \ }
Suppose that $G$ admits a bipartition $(A,B)$ with $|A|\ge |B|>0$ such that every vertex in $A$ has degree at least five and
$G$ has no $K_5$ minor. We may assume that $|A|=|B|$. Choose $b\in B$.
Now take $k$ copies of $G$, say $G_1\LL G_k$, and let $b_i$
be the vertex of $G_i$ that corresponds to $b$. Let $H$ be obtained from the disjoint union of $G_1\LL G_k$ by identifying
$b_1\LL b_k$. Then $H$ has no $K_5$ minor, and has a bipartition $(C,D)$ with $|C|=k|A|$ and $|D|=k(|B|-1)+1=k(|A|-1)+1$,
and every vertex in $C$ has degree at least five.
Now add one more vertex $d$ to $H$ adjacent to every vertex in $C$; then the graph we produce admits a bipartition
$(C,D\cup \{d\})$ where every vertex in $C$ has degree at least six, and it has no $K_6$ minor (because $H$ has no $K_5$
minor). So if we choose $k$ such that $k|A|\ge k(|A|-1)+2$, that is, $k\ge 2$, we obtain a contradiction to \ref{mainthm}.
This proves \ref{smallcases5}.~\bbox

\section{Some lemmas}
Let us begin on the proof of \ref{mainthm}. Thus, let $G$ be a graph that admits a bipartition $(A,B)$ with $|A|\ge |B|>0$
such that every vertex in $A$ has degree at least six; we need to show that $G$ admits a $6$-cluster. 
Let us say $G$ is a {\em candidate} if admits a bipartition $(A,B)$ with $|A|\ge |B|>0$
such that every vertex in $A$ has degree at least six, and $G$ has no 6-cluster. We need to show that there is no candidate. 
We call $(A,B)$ the {\em bipartition} of the candidate.
If $G$ is a candidate with $|G|+|E(G)|$ minimum, we say it is a {\em minimal candidate}.

Let us begin with some easy observations.
\begin{thm}\label{facts}
Let $G$ be a minimal candidate with bipartition $(A,B)$. Then 
\begin{itemize}
\item $|A|=|B|$; 
\item every vertex in $A$ has degree exactly six; 
\item if $X\subseteq B$ is nonempty then $G\setminus X$ has at most $|X|$ components;
\item if $X\subseteq A$ is nonempty then $G\setminus X$ has at most $|X|$ components; and
\item  if $X\subseteq A$ is nonempty and $G\setminus X$ has exactly $|X|$ components then at most one of them has more than one vertex.
\end{itemize}
\end{thm}
\Proof
If $|A|>|B|$ we could delete a vertex in $A$ and obtain a smaller candidate; and if some vertex in $A$ has degree more than six 
we could delete an edge incident with it to obtain a smaller candidate, in either case contradicting minimality.

For the third bullet, let $X\subseteq B$ be nonempty, and let $G_1\LL G_k$ be the components of $G\setminus X$.
For $1\le i\le k$, let $|V(G)\cap A|=p_i$ and $|V(G)\cap B|=q_i$.
Then for $1\le i\le k$, since 
$G\setminus V(G_i)$ is not a candidate, it follows that $|A|-p_i<|B|-q_i$, and so $p_i\ge q_i+1$ since $|A|=|B|$;
but then 
$$|A|=p_1+\cdots+p_k\ge q_1+\cdots+q_k+k=|B|-|X|+k,$$ 
and since
$|A|=|B|$ it follows that $k\le |X|$. This proves the third bullet.

For the fourth bullet, let $X\subseteq A$ be nonempty, and let $G_1\LL G_k$ be the components of $G\setminus X$.
For $1\le i\le k$, let $|V(G)\cap A|=p_i$ and $|V(G)\cap B|=q_i$. For $1\le i\le k$, since $G_i$ is not a candidate, it follows
that $p_i\le q_i-1$; but
$$|A|=p_1+\cdots+p_k+|X|\le q_1+\cdots+q_k-k+|X|=|B|-k+|X|,$$
and since
$|A|=|B|$ it follows that $k\le |X|$. This proves the fourth bullet.

Finally, in the same notation, suppose that $k=|X|$, and $G_1,G_2$ both have at least two vertices. Thus $p_i\le q_i-1$ 
for $1\le i\le k$, and since $k=|X|$, it follows that  $p_i= q_i-1$ for $1\le i\le k$. From the third and fourth 
bullets it follows that $G$ is two-connected, and so there are two vertex-disjoint paths of $G$, say $R,S$ each with first 
vertex in $V(G_1)$ and last vertex in $V(G_2)$, and each with no other vertices in $V(G_1)\cup V(G_2)$. Consequently $R,S$
each have first vertex in $V(G_1)\cap B$ and last vertex in $V(G_2)\cap B$. By contracting $V(R)$ and $V(S)$ into
$B$ we see that $G$ contains as a minor the
graph obtained from $G_1\cup G_2$ by identifying the ends of $R$ and identifying the ends of $S$. But this graph admits a bipartition with parts of cardinalities $p_1+p_2$ and $q_1+q_2-2=p_1+p_2$, and so it is a smaller candidate, a contradiction.
This proves the fifth bullet and so proves \ref{facts}.~\bbox

\ref{facts} has a useful corollary:
\begin{thm}\label{5comps}
Let $G$ be a minimal candidate with bipartition $(A,B)$, and let $X\subseteq A$ or $X\subseteq B$ with $|X|=4$. 
Then there do not exist five 
connected subgraphs $Y_1\LL Y_5$ of $G\setminus X$, pairwise vertex-disjoint, such that for $1\le i\le 5$, every vertex in $X$ 
has a neighbour in $Y_i$.
\end{thm}
\Proof Let $X=\{x_1,x_2,x_3,x_4\}$.
Suppose that such $Y_1\LL Y_5$ exist, and choose them with maximal union. By the third and fourth bullets of \ref{facts} they are not 
all components of $G\setminus X$, and so from the maximality of their union, some two of them are joined by an edge, say $Y_4,Y_5$.
But then there is a 6-cluster 
$$\{V(Y_1)\cup \{x_1\}, V(Y_2)\cup \{x_2\}, V(Y_3)\cup \{x_3\}, V(Y_4), V(Y_5), \{x_4\}\},$$
which is impossible. This proves \ref{5comps}.~\bbox

Here is another way of using the minimality of the candidate. Let $a_1\LL a_p\in A$ be distinct and $b_1\LL b_q\in B$
be distinct. The {\em cover graph} $H$ (with respect to $a_1\LL a_p,b_1\LL b_q$) is the graph with vertex set $\{b_1\LL b_q\}$
in which two distinct vertices $u,v$ are adjacent if there is a vertex $w\in A\setminus \{a_1\LL a_p\}$ adjacent in $G$ 
to both $u,v$. We denote the chromatic number of $H$ by $\chi(H)$. A partition of $V(H)=\{b_1\LL b_q\}$ into sets that are stable in $H$ is a {\em colouring} of $H$, and a partition
$\{Y_1\LL Y_k\}$ of $V(H)$ is 
{\em feasible} if there are pairwise disjoint subsets $X_1\LL X_k$ of $\{a_1\LL a_p\}$ such that $G[X_i\cup Y_i]$ is connected for 
$1\le i\le k$. (Note that the sets $X_i$ might be empty.)
\begin{thm}\label{colouring}
Let $G$ be a minimal candidate with bipartition $(A,B)$, and let $a_1\LL a_p\in A$ be distinct and $b_1\LL b_q\in B$, 
with cover graph $H$. Then no colouring of $H$ of cardinality at most $q-p$ is feasible.
\end{thm}
\Proof
Suppose that the colouring $\{Y_1\LL Y_k\}$ of $H$ is feasible, where $k\le q-p$, and let $X_1\LL X_k$ be the corresponding 
subsets of $\{a_1\LL a_p\}$. By 
contracting each of the sets $X_i\cup Y_i$ into $B$, we obtain a graph with a bipartition $(C,D)$ say, where
$|C|\ge |A|-p$ and $|D|=|B|-q+k\le |B|-p\le |C|$; and every vertex in $C$ has degree at least six, since each of $Y_1\LL Y_k$ 
is stable in $H$, and so this graph is a candidate, which is 
impossible from the minimality of $G$. This proves \ref{colouring}.~\bbox

A special case of \ref{colouring} is used so frequently that it is worth stating explicitly:
\begin{thm}\label{mate}
Let $G$ be a minimal candidate with bipartition $(A,B)$. If $b_1,b_2\in B$ are distinct and have a common neighbour in $A$ then they have at 
least two common neighbours in $A$.
\end{thm}
\Proof
Let $a\in A$ be adjacent to $b_1,b_2$. If $b_1,b_2$ have no other common neighbour, then the covering graph of $a,b_1,b_2$
admits a colouring of cardinality one, which is therefore feasible, contrary to \ref{colouring}.~\bbox

\section{Excluding $K(3,5,0)$- and $K(4,4,1)$-subgraphs}

We will prove a series of results about minimal candidates, which eventually allow to show that there is no such graph. Most of these result
are of the form ``If $G$ is a minimal candidate, then $G$ has no subgraph of the following type'', where the types describe
subgraphs that become smaller and
simpler as the sequence goes on. For instance, one of our result will say that
there do not exist two vertices in $A$ and six vertices
in $B$ such that each of the first is adjacent to each of the second.  We need some notation to describe these ``types''.
For integers $p,q,r\ge 0$ with $r\le \min(p,q)$,
let us say a subgraph $H$
of $G$ is a {\em $K(p,q,r)$-subgraph} if it consists of $p$ vertices $a_1\LL a_p\in A$ and $q$ vertices $b_1\LL b_q\in B$,
where the pairs $a_1b_1, a_2b_2\LL a_rb_r$ are nonadjacent, and otherwise each $a_i$ is adjacent to each $b_j$.
Thus $H$ is obtained from a complete bipartite graph $K_{p,q}$ by deleting a matching with $r$ edges; but it matters that the $p$
vertices belong to $A$ and the $q$ belong to $B$, and not the other way around.

In this section we will prove that a minimal candidate has no $K(3,5,0)$-subgraph and no $K(4,4,1)$-subgraph.
We begin with: 
\begin{thm}\label{440}
Let $G$ be a minimal candidate with bipartition $(A,B)$. Then $G$ has no $K(4,4,0)$-subgraph.
\end{thm}
\Proof
Suppose that $a_1\LL a_4\in A$ are adjacent to $b_1\LL b_4\in B$. Let $Z=\{a_1\LL a_4,b_1\LL b_4\}$. For each component $C$ of 
$G\setminus Z$,
let $N(C)$ denote the set of vertices in $Z$ with a neighbour in $V(C)$.
\\
\\
(1) {\em For each component $C$ of
$G\setminus Z$, $\{a_1,a_2,a_3,a_4\}\not\subseteq N(C)$, and $\{b_1,b_2,b_3,b_4\}\not\subseteq N(C)$.}
\\
\\
This is immediate from two applications of \ref{5comps}, setting $X=\{a_1,a_2,a_3,a_4\}$ and $X=\{b_1,b_2,b_3,b_4\}$.
\begin{figure}[H]
\centering

\begin{tikzpicture}[scale=0.8,auto=left]
\tikzstyle{every node}=[inner sep=1.5pt, fill=black,circle,draw]
\node (a1) at (-3,2) {};
\node (a2) at (-1,2) {};
\node (a3) at (1,2) {};
\node (a4) at (3,2) {};
\node (b1) at (-3,0) {};
\node (b2) at (-1,0) {};
\node (b3) at (1,0) {};
\node (b4) at (3,0) {};

\foreach \from/\to in {a1/b1,a1/b2,a1/b3,a1/b4, a2/b1,a2/b2,a2/b3,a2/b4,a3/b1,a3/b2,a3/b3,a3/b4,a4/b1,a4/b2,a4/b3,a4/b4}
\draw [-] (\from) -- (\to);
\tikzstyle{every node}=[]
\draw[above] (a1) node            {$a_1$};
\draw[above] (a2) node            {$a_2$};
\draw[above] (a3) node            {$a_3$};
\draw[above] (a4) node            {$a_4$};
\draw[below] (b1) node            {$b_1$};
\draw[below] (b2) node            {$b_2$};
\draw[below] (b3) node            {$b_3$};
\draw[below] (b4) node            {$b_4$};

\end{tikzpicture}
\caption{$K(4,4,0)$-subgraph.} \label{fig:440}
\end{figure}
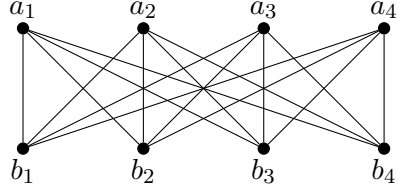

\bigskip

Since $a_1,a_2,a_3,a_4$ have degree six, and so each of them belongs to $N(C)$ for some component $C$ of $G\setminus Z$, it
follows from (1) that there are at least two such components.
\\
\\
(2) {\em For each component $C$ of
$G\setminus Z$, $|N(C)\cap A|\ne 1$.}
\\
\\
Suppose that $N(C)\cap A=\{a_1\}$ say. Let $b\in V(C)$ be adjacent to $a_1$; thus $b\in B$. By (1) we may assume that 
$b_1\notin N(C)$. But then $b,b_1$ have a unique common neighbour, contrary to \ref{mate}. This proves (2).
\\
\\
(3) {\em For each component $C$ of
$G\setminus Z$, one of $|N(C)\cap A|, |N(C)\cap B|\ge 2$.}
\\
\\
If $|N(C)\cap A|\le 1$ then $N(C)\cap A=\emptyset$ by (2), and so by the third bullet of \ref{facts}, taking $X=N(C)$,
it follows that $|N(C)\cap B|\ge 2$. This proves (3).
\\
\\
(4) {\em $|N(C)\cap B|\le 1$ and $|N(C)\cap A|\in \{2,3\}$ for every component $C$ of $G\setminus Z$.}
\\
\\
Suppose that $|N(C)\cap B|\ge 2$, and let $b_1,b_2\in N(C)$ say. By (1), there is a component $C'\ne C$ of $G\setminus Z$
with $N(C')\cap A\ne \emptyset$, and hence with $|N(C')\cap A|\ge 2$ by (2). Let $a_1,a_2\in N(C')$ say.
Then there is a 6-cluster
$$\{\{a_1\},V(C')\cup\{a_2\},\{b_1\},V(C)\cup \{b_2\},\{a_3,b_3\}, \{a_4,b_4\}\},$$
a contradiction. This proves that $|N(C)\cap B|\le 1$ for every component $C$ of $G\setminus Z$; and so 
$|N(C)\cap A|\in \{2,3\}$ for every component $C$ of $G\setminus Z$ by (3) and (1). This proves (4).
\\
\\
(5) {\em If $v\in B\setminus Z$ has a neighbour in $\{a_1,a_2,a_3,a_4\}$ then it has at least two such neighbours.}
\\
\\
Let $C$ be the component of $G\setminus Z$ that contains $v$. We may assume that $v$ is adjacent to $a_1$, and by 
(1), we may assume that $b_1\notin N(C)$. By \ref{mate}, $b_1,v$ have another common neighbour, which must be in $\{a_2,a_3,a_4\}$
since $b_1\notin N(C)$. This proves (5).
\\
\\
(6) {\em If $u,v\in B\setminus Z$ have a common neighbour in $\{a_1,a_2,a_3,a_4\}$ and belong to different components of
$G\setminus Z$ then they have at least two common neighbours in $\{a_1,a_2,a_3,a_4\}$. Consequently, if $C,C'$
are distinct components of $G\setminus Z$ then $|N(C)\cap N(C')\cap A|\ne 1$.}
\\
\\
The first claim follows from \ref{mate}  applied to $u,v$; and the second is a consequence. This proves~(6).
\\
\\
(7) {\em $|N(C)\cap A|=2$ for each component $C$ of $G\setminus Z$.}
\\
\\
Suppose not; then by (4) $|N(C)\cap A|=3$, and we may assume that $N(C)\cap A=\{a_1,a_2,a_3\}$. Let $C'$ be a component
of $G\setminus Z$ with $a_4\in N(C')$. By (4), $N(C)\cap N(C')\ne \emptyset$, and so by (6), $|N(C')\cap A|= 3$ and we 
may assume that $a_2,a_3,a_4\in N(C')$. Since $a_2$ has a neighbour in each of $C,C'$ and is also adjacent to $b_1,b_2,b_3,b_4$,
it has no more neighbours, and the same holds for $a_3$. Consequently if $C''\ne C,C'$ is a component of $G\setminus Z$
then $N(C'')\cap A\subseteq \{a_1,a_4\}$, and so equality holds by (4), contrary to (6). Thus $C,C'$ are the only components
of $G\setminus Z$.
Hence $a_1$ has two neighbours $d_1,d_2\in V(C)$, and $a_2,a_3$ each have exactly one neighbour in $V(C')$. By (5), 
each of $d_1,d_2$ is adjacent to two of $a_1,a_2,a_3,a_4$, and so we may assume that $d_1$ is adjacent to $a_2$ and not to $a_3$.
Similarly there exists $d'\in V(C')$ adjacent to $a_2$ and not to $a_3$. But then $d_1,d'$ have a unique common neighbour,
contrary to (6). This proves (7).

\bigskip

By (4), $b_1,b_2,b_3,b_4$ belong to different components of $G\setminus \{a_1,a_2,a_3,a_4\}$, and so by \ref{facts},
these are the only components of $G\setminus \{a_1,a_2,a_3,a_4\}$, and three of them have only one vertex. Consequently
we may assume that $b_2,b_3,b_4$ have degree four in $G$, and $b_1\in N(C)$ for each component $C$ of $G\setminus Z$.
By (7) and (6), we may assume that for each component $C$ of $G\setminus Z$, $N(C)\cap A=\{a_1,a_2\}$ or $\{a_3,a_4\}$.
Let $G_1$ be the union of the components $C$ with $N(C)\cap A=\{a_1,a_2\}$, and define $G_2$ similarly for $\{a_3,a_4\}$.
Thus $a_1,a_2$ each have two neighbours in $V(G_1)$, and $a_3,a_4$ each have two in $V(G_2)$. Hence by contracting 
 $\{a_1,b2,a_3\}$ and $\{a_2,b_3, b_4\}$ into $A$ we obtain a smaller candidate, a contradiction.
This proves \ref{440}.~\bbox

\begin{thm}\label{360}
Let $G$ be a minimal candidate with bipartition $(A,B)$. Then $G$ has no $K(3,6,0)$-subgraph.
\end{thm}
\Proof
Suppose that $a_1,a_2,a_3\in A$ are all adjacent to each of $b_1\LL b_6\in B$. Let $H$ be the cover graph with respect to 
 $a_1,a_2,a_3, b_1\LL b_6$. 
\begin{figure}[H]
\centering

\begin{tikzpicture}[scale=0.8,auto=left]
\tikzstyle{every node}=[inner sep=1.5pt, fill=black,circle,draw]
\node (a1) at (-3,3) {};
\node (a2) at (-0,3) {};
\node (a3) at (3,3) {};
\node (b1) at (-5,0) {};
\node (b2) at (-3,0) {};
\node (b3) at (-1,0) {};
\node (b4) at (1,0) {};
\node (b5) at (3,0) {};
\node (b6) at (5,0) {};

\foreach \from/\to in {a1/b1,a1/b2,a1/b3,a1/b4, a1/b5,a1/b6, a2/b1,a2/b2,a2/b3,a2/b4,a2/b5,a2/b6,a3/b1,a3/b2,a3/b3,a3/b4,a3/b5,a3/b6}
\draw [-] (\from) -- (\to);
\tikzstyle{every node}=[]
\draw[above] (a1) node            {$a_1$};
\draw[above] (a2) node            {$a_2$};
\draw[above] (a3) node            {$a_3$};

\draw[below] (b1) node            {$b_1$};
\draw[below] (b2) node            {$b_2$};
\draw[below] (b3) node            {$b_3$};
\draw[below] (b4) node            {$b_4$};
\draw[below] (b5) node            {$b_5$};
\draw[below] (b6) node            {$b_6$};

\end{tikzpicture}
\caption{$K(3,6,0)$-subgraph.} \label{fig:360}
\end{figure}
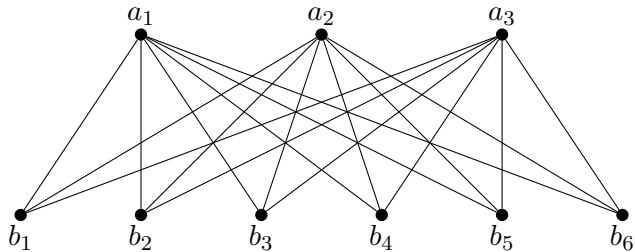
\noindent
(1) {\em If $b_1,b_2,b_3$ are pairwise adjacent in $H$, there is a vertex $a\ne a_1,a_2,a_3$ adjacent to all of $b_1,b_2,b_3$,
and no other vertex in $A\setminus \{a_1,a_2,a_3\}$ is adjacent to any two of $b_1,b_2,b_3$.}
\\
\\
Since $b_1b_2\in E(H)$, there exists $c_3\in A\setminus \{a_1,a_2,a_3\}$ adjacent to $b_1,b_2$; and similarly there exists $c_1$
adjacent to $b_2,b_3$ and $c_2$ adjacent to $b_3,b_1$. If $c_1,c_2,c_3$ are all different, there is a 6-cluster
$$\{\{c_2,b_1\},\{c_3,b_2\},\{c_1,b_3\},\{a_1,b_4\},\{a_2,b_5\},\{a_3,b_6\}\},$$
a contradiction. So $c_1,c_2,c_3$ cannot be chosen all different. In particular we may assume that some $c\in A\setminus \{a_1,a_2,a_3\}$ is adjacent to 
$b_1,b_2,b_3$. If some other vertex $d \in A\setminus \{a_1,a_2,a_3\}$ is adjacent to two of $b_1,b_2,b_3$, say to $b_1,b_2$,
 there is a 6-cluster
$$\{\{b_1,d\},\{b_2\},\{c,b_3\},\{a_1,b_4\},\{a_2,b_5\},\{a_3,b_6\}\},$$
a contradiction. This proves (1).

\bigskip

Every colouring of $H$ of cardinality three is feasible (as we can add one $a_i$ to each vertex class), so $\chi(H)\ge 4$ by \ref{colouring}. Consequently either $H$ consists 
of an induced cycle of length five together with one more vertex adjacent to every vertex of the cycle, or $H$ has a clique of 
size four. In either case there are four vertices of $H$ such that five of the six pairs of them are adjacent in $H$.
We may assume that $b_1b_2,b_1b_3,b_1b_4,b_2b_3,b_2b_4$ are all edges of $H$. By (1) there exists $c\in A\setminus \{a_1,a_2,a_3\}$ adjacent to
 $b_1,b_2,b_3$, and $d\in A\setminus \{a_1,a_2,a_3\}$ adjacent to
 $b_1,b_2,b_4$; and by (1) again, $c=d$. Thus $c$ is adjacent to $b_1,b_2,b_3,b_4$, and so $G[\{a_1,a_2,a_3,c,b_1,b_2,b_3,b_4\}]$
is a $K(4,4,0)$-subgraph, contrary to \ref{440}. This proves \ref{360}.~\bbox

If $P$ is a path, we denote by $P^*$ the set of vertices in the interior of $P$, that is, the vertices that have degree two in $P$.
\begin{thm}\label{350}
Let $G$ be a minimal candidate with bipartition $(A,B)$. Then $G$ has no $K(3,5,0)$-subgraph.
\end{thm}
\Proof
Suppose that $a_1,a_2,a_3\in A$ are all adjacent to each of $b_1\LL b_5\in B$. 
\begin{figure}[H]
\centering

\begin{tikzpicture}[scale=0.8,auto=left]
\tikzstyle{every node}=[inner sep=1.5pt, fill=black,circle,draw]
\node (a1) at (-3,3) {};
\node (a2) at (-0,3) {};
\node (a3) at (3,3) {};
\node (b1) at (-4,0) {};
\node (b2) at (-2,0) {};
\node (b3) at (-0,0) {};
\node (b4) at (2,0) {};
\node (b5) at (4,0) {};

\foreach \from/\to in {a1/b1,a1/b2,a1/b3,a1/b4, a1/b5, a2/b1,a2/b2,a2/b3,a2/b4,a2/b5,a3/b1,a3/b2,a3/b3,a3/b4,a3/b5}
\draw [-] (\from) -- (\to);
\tikzstyle{every node}=[]
\draw[above] (a1) node            {$a_1$};
\draw[above] (a2) node            {$a_2$};
\draw[above] (a3) node            {$a_3$};

\draw[below] (b1) node            {$b_1$};
\draw[below] (b2) node            {$b_2$};
\draw[below] (b3) node            {$b_3$};
\draw[below] (b4) node            {$b_4$};
\draw[below] (b5) node            {$b_5$};

\end{tikzpicture}
\caption{$K(3,5,0)$-subgraph.} \label{fig:350}
\end{figure}
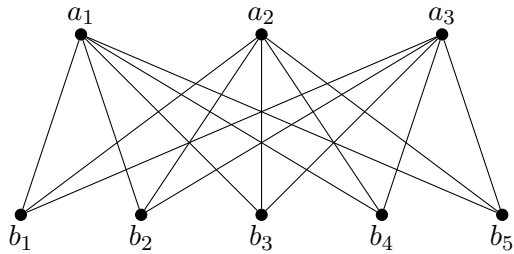
Each of $a_1,a_2,a_3$ has exactly one neighbour different from $b_1\LL b_5$, and they are not all equal since $G$ has no
$K(3,6,0)$-subgraph by \ref{360}. So some vertex is adjacent to exactly one of $a_1,a_2,a_3$; say $b_6$ is adjacent to $a_1$
and not to $a_2,a_3$. By \ref{mate}, for $1\le i\le 5$ $b_i,b_6$ have a common neighbour different from $a_1$. Choose a set $X$
of neighbours of $b_6$, with $a_1,a_2,a_3\in X$, minimal such that $b_1\LL b_5$ each have a neighbour in $X$. 
Consequently for each $x\in X$ there exists $i\in \{1\LL 5\}$ such that $x$ is the unique neighbour of $b_i$ in $X$.
\\
\\
(1) {\em Every vertex different from $a_1,a_2,a_3$ with two neighbours in $\{b_1\LL b_5\}$ is in $X$.}
\\
\\
Suppose that $a\in A\setminus \{a_1,a_2,a_3\}$ is adjacent to $a_1,a_2$ say, and $a\notin X$. Then there is a 6-cluster
$$\{\{b_1\},\{a,b_2\}, \{b_6,b_3\}\cup X,\{a_1\},\{a_2,b_4\},\{a_3,b_5\}\},$$
a contradiction. This proves (1).
\\
\\
(2) {\em Some vertex different from $a_1,a_2,a_3$ has three neighbours in $\{b_1\LL b_5\}$.}
\\
\\
Let $H$ be the cover graph with respect to
$a_1,a_2,a_3, b_1\LL b_5$. By \ref{colouring}, $\chi(H)\ge 3$. So either it is a cycle of
length five, or it has a triangle.  (A {\em triangle} means a clique with
cardinality three.) Suppose first that $H$ is a cycle of length five, with edges $b_1b_2,b_2b_3,b_3b_4,b_4b_5,b_5b_1$
say. Some vertex $d_{1,2}\ne a_1,a_2,a_3$ is adjacent in $G$ to $b_1,b_2$, from the definition of $H$, and it is 
nonadjacent to $b_3,b_4,b_5$ since $H$ is a cycle. Define $d_{2,3}$ and so on similarly. By (1), each of these five vertices
is in the set $X$; but then none of $b_1\LL b_5$ has a unique neighbour in $X$, contrary to the minimality of $X$.

It follows that $H$ has a triangle, say with vertices $b_1,b_2,b_3$. 
Some vertex $d_{1,2}\ne a_1,a_2,a_3$ is adjacent in $G$ to $b_1,b_2$, from the definition of $H$; define $d_{2,3},d_{3,1}$
similarly. Suppose that $d_{1,2}, d_{2,3}, d_{3,1}$ are all different. Then there is a 6-cluster
$$\{\{d_{1,2},b_1\},\{d_{2,3},b_2\},\{d_{3,1},b_3\},\{a_1\}, \{a_2,b_4\}, \{a_3,b_5\}\},$$
a contradiction. So two of $d_{1,2}, d_{2,3}, d_{3,1}$ are equal. This proves (2).

\bigskip

Let $a_4$ be adjacent to $b_1,b_2,b_3$ say. It is nonadjacent to $b_4,b_5$ since $G$ has no $K(4,4,0)$-subgraph by \ref{440}.
By \ref{colouring} the cover graph with respect to $a_1,a_2,a_3, a_4, b_1\LL b_5$ has chromatic number at least two, and so 
has an edge. Choose $a_5$ different from $a_1\LL a_4$ with two neighbours in $\{b_1\LL b_5\}$. By (1), $a_4,a_5$
are both adjacent to $b_6$.
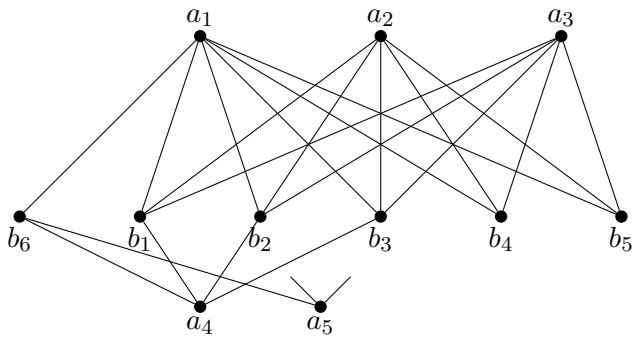
\begin{figure}[H]
\centering

\begin{tikzpicture}[scale=0.8,auto=left]
\tikzstyle{every node}=[inner sep=1.5pt, fill=black,circle,draw]
\node (a1) at (-3,3) {};
\node (a2) at (-0,3) {};
\node (a3) at (3,3) {};
\node (b1) at (-4,0) {};
\node (b2) at (-2,0) {};
\node (b3) at (-0,0) {};
\node (b4) at (2,0) {};
\node (b5) at (4,0) {};
\node (b6) at (-6,0) {};
\node (a4) at (-3,-1.5) {};
\node (a5) at (-1,-1.5) {};

\foreach \from/\to in {a1/b1,a1/b2,a1/b3,a1/b4, a1/b5, a2/b1,a2/b2,a2/b3,a2/b4,a2/b5,a3/b1,a3/b2,a3/b3,a3/b4,a3/b5, a1/b6, a4/b1,
a4/b2,a4/b3,a4/b6, a5/b6}
\draw [-] (\from) -- (\to);
\draw[-] (a5) to (-1.5,-1);
\draw[-] (a5) to (-.5,-1);
\tikzstyle{every node}=[]
\draw[above] (a1) node            {$a_1$};
\draw[above] (a2) node            {$a_2$};
\draw[above] (a3) node            {$a_3$};

\draw[below] (b1) node            {$b_1$};
\draw[below] (b2) node            {$b_2$};
\draw[below] (b3) node            {$b_3$};
\draw[below] (b4) node            {$b_4$};
\draw[below] (b5) node            {$b_5$};
\draw[below] (b6) node            {$b_6$};
\draw[below] (a4) node            {$a_4$};
\draw[below] (a5) node            {$a_5$};

\end{tikzpicture}
\caption{For the last part of the proof of \ref{350}. $a_5$ is adjacent to two of $b_1\LL b_5$.} \label{fig:350+}
\end{figure}

Up to symmetry there are three cases: $a_5$ is adjacent to $b_1,b_2$; $a_5$ is adjacent to $b_1,b_4$; and $a_5$ is adjacent to
$b_4,b_5$. 

First, if $a_5$ is adjacent to $b_1,b_2$, there is a 6-cluster
$$\{\{b_1\},\{a_5,b_2\},\{a_4,b_3\},\{a_1\},\{a_2,b_4\},\{a_3,b_5\}\},$$
a contradiction. If $a_5$ is adjacent to $b_1,b_4$, there is a 6-cluster
$$\{\{b_1\},\{a_4,b_2\},\{a_5,b_4,b_6\},\{a_1\},\{a_2,b_3\},\{a_3,b_5\}\},$$
a contradiction. So $a_5$ is adjacent to $b_4,b_5$. 

By \ref{facts}, the graph $G\setminus \{a_1,a_2,a_3,a_4,a_5\}$ has at most five components; and so some two of $b_1\LL b_6$
belong to the same component. So there is a path $P$ of $G$ between two of  $b_1\LL b_6$ with no other vertices
in $\{a_1\LL b_5,b_1\LL b_6\}$. Let $P$ have ends $b_i, b_j$ say. The subgraph induced on $\{a_4,a_5,b_1\LL b_6\}$ is a tree,
and its union with $P$ includes a cycle that contains $P$; and in all cases we can use this cycle to make
three of $b_1\LL b_5$ adjacent and thereby produce a $K_6$ minor. 
In detail (up to symmetry these are the only possibilities):
\begin{itemize}
\item If $(i,j)=(1,2)$, there is a 6-cluster
$\{\{b_1\},\{P^*\cup \{b_2\}\},\{a_4,b_3\},\{a_1\},\{a_2,b_4\},\{a_3,b_5\}\}.$
\item If $(i,j)=(1,4)$, there is a 6-cluster
$\{\{b_1\},\{a_4,b_2\},P^*\cup \{b_4,b_6\},\{a_1\},\{a_2,b_3\},\{a_3,b_5\}\}.$
\item If $(i,j)=(1,6)$ there is a 6-cluster
$\{\{b_1\},\{a_4b_2\}, P^*\cup \{b_6,a_5,b_4\},\{a_1\},\{a_2,b_3\},\{a_3,b_5\}\}.$
\item If $(i,j) = (4,5)$ there is a 6-cluster
$\{\{b_1,a_4,b_6,a_5\},\{b_4\},P^*\cup \{b_5\},\{a_1\},\{a_2,b_2\},\{a_3,b_3\}\}.$
\item If $(i,j) = (4,6)$ there is a 6-cluster
$\{P^*\cup \{b_6\},\{b_4\},\{a_5,b_5\},\{a_1\},\{a_2,b_2\},\{a_3,b_3\}\}.$
\end{itemize}
In each case we have a contradiction. This proves \ref{350}.~\bbox

\begin{thm}\label{474}
Let $G$ be a minimal candidate with bipartition $(A,B)$. Then $G$ has no $K(4,7,4)$-subgraph.
\end{thm}
\Proof
Suppose that $a_1\LL a_4\in A$ are all adjacent to each of $b_1\LL b_7\in B$, except the pairs $a_1b_1,a_2b_2,a_3b_3,a_4b_4$.
Let $H$ be the cover graph with respect to $a_1\LL a_4,b_1\LL b_7$. We claim that every partition of $V(H)$ into at most three 
sets is 
feasible. To see this, let $\{Y_1\LL Y_k\}$ be a partition of $V(H)$ with $k\le 3$. 
We must show there are disjoint subsets
$X_1\LL X_k$ of $\{a_1\LL a_4\}$ such that $X_i\cup Y_i$ is connected for $1\le i\le k$.
If some $Y_i$, say $Y_1$,
contains all of $b_1\LL b_4$ we may set
$X_1=\{a_1,a_2\}$ and $X_2\LL X_k$ each to contain one of $a_3,a_4$; so we may assume that 
for $1\le i\le k$, there exists $j\in \{1\LL 4\}$ such that $b_j\notin Y_i$. But then (from Hall's ``marriage'' theorem, for instance),
there is an injection $\phi:\{1\LL k\}\rightarrow\{1\LL 4\}$ such that $b_{\phi(i)}\notin Y_i$ for $1\le i\le k$; so we may
set $X_i=\{a_{\phi(i)}\}$ for $1\le i\le k$.
This proves that  every partition of $V(H)$ into at most three
sets is
feasible, and so 
$\chi(H)\ge 4$ by \ref{colouring}. 
\\
\\
(1) {\em $b_1,b_2,b_3,b_4$ are pairwise adjacent in $H$, and $H$ has no other edges.}
\\
\\
Suppose that say $b_5b_6$ are adjacent in $H$, and let $P$ be a path of $G$ between $b_5,b_6$ with no other vertex in 
$\{a_1\LL a_4,b_1\LL b_7\}$. There is a 6-cluster
$$\{a_1,b_2\},\{a_2,b_3\},\{a_3,b_4\},\{a_4,b_1\},\{b_5\},P^*\cup \{b_6\}\},$$
a contradiction. So $b_5,b_6,b_7$ are pairwise nonadjacent in $H$. 
Let $S$ be the set of edges of $H$ with both ends in $\{b_1\LL b_4\}$, and $T$ the set of edges of $H$ with one end in
$\{b_1\LL b_4\}$ and the other in $\{b_5,b_6,b_7\}$. Thus every edge of $H$ belongs to exactly one of $S,T$. 
If say $b_3b_4$ and  $b_4b_5$ are edges of $H$,
let $P,Q$ be the corresponding paths of $G$; then there is a 6-cluster
$$\{\{a_1,b_2\},\{a_2,b_1\},\{a_3,b_6\},\{a_4,b_7\},V(P), Q^*\cup \{b_5\}\},$$
a contradiction. Hence no edge in $S$ shares an end with an edge in $T$, and so no component of $H$
has an edge in $S$ and an edge in $T$. But some component $H'$ of $H$ has chromatic number at least four, and so
not all its edges are in $T$, since then it would be bipartite. Hence all its edges are in $S$, and so $V(H')\subseteq \{b_1\LL b_4\}$;
and since $H'$ has chromatic number four, $H'$ is a complete graph. This proves (1).

\begin{figure}[H]
\centering

\begin{tikzpicture}[scale=0.8,auto=left]
\tikzstyle{every node}=[inner sep=1.5pt, fill=black,circle,draw]
\node (a1) at (-4,3) {};
\node (a2) at (-2,3) {};
\node (a3) at (0,3) {};
\node (a4) at (2,3) {};
\node (b1) at (-4,0) {};
\node (b2) at (-2,0) {};
\node (b3) at (-0,0) {};
\node (b4) at (2,0) {};
\node (b5) at (-4,6) {};
\node (b6) at (-1,6) {};
\node (b7) at (2,6) {};

\foreach \from/\to in {a1/b2,a1/b3,a1/b4,a1/b5,a1/b6,a1/b7, a2/b1,a2/b3,a2/b4,a2/b5,a2/b6,a2/b7,a3/b1,a3/b2,a3/b4,
a3/b5, a3/b6,a3/b7,a4/b1,a4/b2,a4/b3,a4/b5,a4/b6,a4/b7}
\draw [-] (\from) -- (\to);
\tikzstyle{every node}=[]
\draw[left] (a1) node            {$a_1$};
\draw[left] (a2) node            {$a_2$};
\draw[right] (a3) node            {$a_3$};
\draw[right] (a4) node            {$a_4$};

\draw[below] (b1) node            {$b_1$};
\draw[below] (b2) node            {$b_2$};
\draw[below] (b3) node            {$b_3$};
\draw[below] (b4) node            {$b_4$};
\draw[above] (b5) node            {$b_5$};
\draw[above] (b6) node            {$b_6$};
\draw[above] (b7) node            {$b_7$};

\end{tikzpicture}
\caption{$K(4,7,4)$-subgraph.} \label{fig:474}
\end{figure}
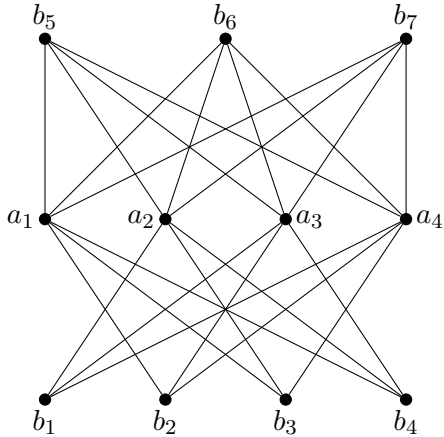

Choose $c\ne a_1\LL a_4$ adjacent to $b_1,b_2$, and $c'\ne a_1\LL a_4$ adjacent to $b_3,b_4$. Then 
$c=c'$, since otherwise the connected subgraphs with vertex sets $\{b_1,c,b_2\},\{b_3,c',b_4\},\{b_5\},\{b_6\},\{b_7\}$
violate \ref{5comps}. Consequently $c$ is adjacent to $b_1,b_2,b_3,b_4$. By the same argument, no other vertex has two neighbours
in $\{b_1\LL b_4\}$, and hence  no other vertex has two neighbours
in $\{b_1\LL b_7\}$, since $b_5,b_6,b_7$ have degree zero in $H$. But then the cover graph with respect to 
$a_1\LL a_4, c, b_1\LL b_7$ has no edges and \ref{colouring} is violated. This proves \ref{474}.~\bbox

\begin{thm}\label{452}
Let $G$ be a minimal candidate with bipartition $(A,B)$. Then $G$ has no $K(4,5,2)$-subgraph.
\end{thm}
\Proof
Suppose that $a_1\LL a_4\in A$ are all adjacent to each of $b_1\LL b_5\in B$, except the pairs $a_1b_1,a_2b_2$.
\\
\\
(1) {\em Every vertex in $B\setminus \{b_1\LL b_5\}$ with a neighbour in $\{a_1\LL a_4\}$ has exactly two neighbours in this set.}
\\
\\
Suppose that some vertex $b_6\in B\setminus \{b_1\LL b_5\}$ is adjacent to exactly one of $a_1\LL a_4$. By \ref{mate}, 
for $i = 3,4,5$ there is a vertex
$c_i$ adjacent to $b_i,b_6$ and not in $\{a_1\LL a_4\}$. But then
$$\{\{a_1,b_2\},\{a_2,b_3\},\{a_3,b_1\},\{b_4\},\{a_4\},\{b_5,c_5,b_6,c_4\}\}$$
is a 6-cluster, a contradiction.

So every vertex in $B\setminus \{b_1\LL b_5\}$ with a neighbour in $\{a_1\LL a_4\}$ has at least two neighbours in this set.
No vertex is adjacent to all of $a_1,a_3,a_4$ or to all of $a_2,a_3,a_4$, since there is no $K(3,5,0)$-subgraph by \ref{350}.
If some vertex $b_6$ different from $b_1\LL b_5$ is adjacent to $a_1,a_2,a_3$, then each of $a_1,a_2,a_4$ has exactly one 
neighbour different
from $b_1\LL b_6$, and these must all be equal since no vertex has one neighbour in $\{a_1,a_2,a_3,a_4\}$; but then $G$
has a $K(4,7,4)$-subgraph, contrary to \ref{474}. Similarly no vertex different from $b_1\LL b_5$ is adjacent to $a_1,a_2,a_4$;
and so  every vertex in $B\setminus \{b_1\LL b_5\}$ with a neighbour in $\{a_1\LL a_4\}$ has exactly two neighbours in this set.
This proves (1).
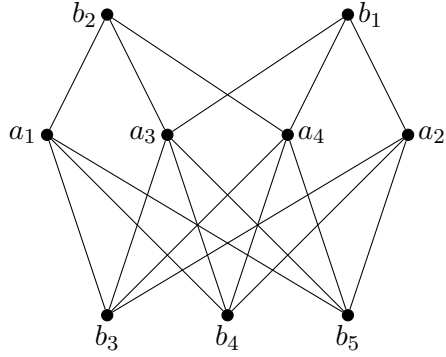
\begin{figure}[H]
\centering

\begin{tikzpicture}[scale=0.8,auto=left]
\tikzstyle{every node}=[inner sep=1.5pt, fill=black,circle,draw]
\node (a2) at (2,3) {};
\node (a1) at (-4,3) {};
\node (a3) at (-2,3) {};
\node (a4) at (0,3) {};
\node (b2) at (-3,5) {};
\node (b1) at (1,5) {};
\node (b3) at (-3,0) {};
\node (b4) at (-1,0) {};
\node (b5) at (1,0) {};

\foreach \from/\to in {a1/b2,a1/b3,a1/b4,a1/b5, a2/b1,a2/b3,a2/b4,a2/b5,a3/b1,a3/b2,a3/b3,a3/b4,
a3/b5, a4/b1,a4/b2,a4/b3,a4/b4,a4/b5}
\draw [-] (\from) -- (\to);
\tikzstyle{every node}=[]
\draw[left] (a1) node            {$a_1$};
\draw[right] (a2) node            {$a_2$};
\draw[left] (a3) node            {$a_3$};
\draw[right] (a4) node            {$a_4$};

\draw[right] (b1) node            {$b_1$};
\draw[left] (b2) node            {$b_2$};
\draw[below] (b3) node            {$b_3$};
\draw[below] (b4) node            {$b_4$};
\draw[below] (b5) node            {$b_5$};

\end{tikzpicture}
\caption{$K(4,5,2)$-subgraph.} \label{fig:452}
\end{figure}
\noindent
(2) {\em $b_3,b_4,b_5$ each have degree four in $G$.}
\\
\\
There are four edges between $\{a_1,a_2\}$ and $B\setminus \{b_1\LL b_5\}$, and only two between 
$\{a_3,a_4\}$ and $B\setminus \{b_1\LL b_5\}$; and so by (1) there is a vertex $b_6\in B\setminus \{b_1\LL b_5\}$ adjacent to $a_1,a_2$
and not to $a_3,a_4$. By \ref{mate}, there is a vertex different from $a_1\LL a_4$ adjacent to $b_6,b_1$ (because
$b_6,b_1$ are adjacent to $a_2$), and similarly there is a vertex different from $a_1\LL a_4$ adjacent to $b_6,b_2$. Consequently
there is a path $P$ between $b_1,b_2$ with no other vertices in $\{a_1\LL a_4,b_1\LL b_5\}$ (possibly containing $b_6$).

We claim that $P,b_3,b_4,b_5$ all belong to different components of $G\setminus \{a_1,a_2,a_3,a_4\}$. Because suppose not;
then, either there is a path $Q$ of $G$ between two of $b_3,b_4,b_5$ with no other vertex in $\{a_1\LL a_4,b_1\LL b_5\}$, or
there is a path $Q$ between one of $b_1,b_2$ and one of $b_3,b_4,b_5$ with no other vertex in $\{a_1\LL a_4,b_1\LL b_5\}$.
In the first case,
say $Q$ has ends $b_3,b_4$; then 
$$\{\{b_3\},Q^*\cup \{b_4\},\{a_4\},\{a_1,b_2\},\{a_2,b_5\},\{a_3,b_1\}\}$$
is a 6-cluster, a contradiction. In the second case, let $Q$ have ends $b_1,b_3$ say; then 
$$\{\{P^*\cup Q^*\cup \{b_1\},\{a_1,b_2\},\{b_3\},\{a_2,b_4\},\{a_3,b_5\},\{a_4\}\}$$
is a 6-cluster, a contradiction. 

This proves that $P,b_3,b_4,b_5$ all belong to different components of $G\setminus \{a_1,a_2,a_3,a_4\}$. By \ref{5comps},
the components containing $b_3,b_4,b_5$ are all singletons. This proves (2).

\bigskip

Now for $i = 3,4$, $a_i$ has one neighbour not in $\{b_1\LL b_5\}$, say $c_i$. Either $c_3=c_4$ and $c_3$
has no neighbour in $\{a_1,a_2\}$, or each of $c_3,c_4$ has a unique neighbour in $\{a_1,a_2\}$, and not the same one.
Thus in either case we may assume that $c_3$ is not adjacent to $a_1$ and $c_4$ is not adjacent to $a_2$. Consequently,
there are at least six vertices in $B\setminus \{b_1,b_2\}$ with a neighbour in $\{a_1,a_3\}$, namely
$b_3, b_4, b_5, c_3$ and the two neighbours of $a_1$ not in $\{b_1\LL b_5\}$. Similarly 
there are at least six vertices in $B\setminus \{b_1,b_2\}$ with a neighbour in $\{a_2,a_4\}$. Since $b_1,b_2$ both have degree
four, by contracting $\{a_1,b_1,a_3\}$ and $\{a_2,b_2,a_4\}$ into $A$ we obtain a smaller candidate, a contradiction.
This proves \ref{452}.~\bbox

\begin{thm}\label{441}
Let $G$ be a minimal candidate with bipartition $(A,B)$. Then $G$ has no $K(4,4,1)$-subgraph.
\end{thm}
\Proof
Suppose that $a_1\LL a_4\in A$ are all adjacent to each of $b_1\LL b_4\in B$, except the pair $a_1b_1$.
No vertex in $B\setminus \{b_1\LL b_4\}$ is adjacent to all of $a_2,a_3,a_4$ since $G$ has no $K(3,5,0)$-subgraph by
\ref{350}. No vertex is adjacent to $a_1$ and to two of $a_2,a_3,a_4$ since  $G$ has no $K(4,5,2)$-subgraph by
\ref{452}. So 
every vertex in $B\setminus \{b_1\LL b_4\}$ with a neighbour in $\{a_1\LL a_4\}$ has at most two neighbours in this set.

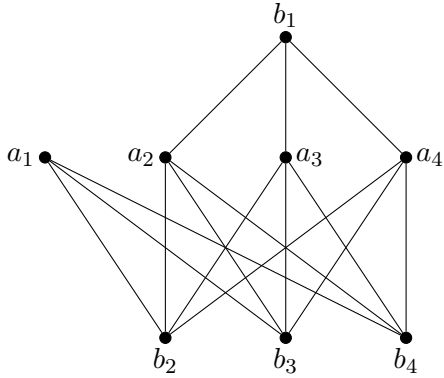
\begin{figure}[H]
\centering

\begin{tikzpicture}[scale=0.8,auto=left]
\tikzstyle{every node}=[inner sep=1.5pt, fill=black,circle,draw]
\node (a1) at (-4,3) {};
\node (a2) at (-2,3) {};
\node (a3) at (0,3) {};
\node (a4) at (2,3) {};
\node (b1) at (0,5) {};
\node (b2) at (-2,0) {};
\node (b3) at (0,0) {};
\node (b4) at (2,0) {};

\foreach \from/\to in {a1/b2,a1/b3,a1/b4,a2/b1,a2/b2,a2/b3,a2/b4,a3/b1,a3/b2,a3/b3,a3/b4,
 a4/b1,a4/b2,a4/b3,a4/b4}
\draw [-] (\from) -- (\to);
\tikzstyle{every node}=[]
\draw[left] (a1) node            {$a_1$};
\draw[left] (a2) node            {$a_2$};
\draw[right] (a3) node            {$a_3$};
\draw[right] (a4) node            {$a_4$};

\draw[above] (b1) node            {$b_1$};
\draw[below] (b2) node            {$b_2$};
\draw[below] (b3) node            {$b_3$};
\draw[below] (b4) node            {$b_4$};

\end{tikzpicture}
\caption{$K(4,4,1)$-subgraph.} \label{fig:441}
\end{figure}
\noindent
Suppose that some  vertex $b_5\in B\setminus \{b_1\LL b_4\}$ has only one neighbour in $\{a_1\LL a_4\}$, and that neighbour is 
different from $a_1$; let it be $a_2$ say.
By \ref{mate}, for $i = 1,2,3,4$ there is a vertex $c_i\ne a_2$ adjacent to both $b_5,b_i$.
But then 
$$\{\{b_5,c_1,c_2,c_3,c_4\}, \{a_2\},\{a_3,b_1\},\{a_4,b_2\},\{a_1,b_3\},\{b_4\}\}$$
is a 6-cluster. So every vertex in $B\setminus \{b_1\LL b_4\}$ with a neighbour in $\{a_2\LL a_4\}$ has exactly two
neighbours in $\{a_1\LL a_4\}$. 

But there are an odd number of edges (nine) between $\{a_1\LL a_4\}$ and $B\setminus \{b_1\LL b_4\}$;
so some vertex $b_5\in B\setminus \{b_1\LL b_4\}$ has a unique neighbour in $\{a_1\LL a_4\}$, and consequently this neighbour is $a_1$.
It follows that 
at most two neighbours of $a_1$ are different from $b_2, b_3, b_4$ and have a neighbour in $\{a_2,a_3,a_4\}$.
But there are six edges between $\{a_2,a_3,a_4\}$ and $B\setminus \{b_1\LL b_4\}$,
and so there is a vertex 
$b_6\in B\setminus \{b_1\LL b_4\}$ adjacent to two of $a_2,a_3,a_4$, say $a_2,a_3$. 
By \ref{mate}, 
for $i = 2,3,4$ there is a vertex $c_i\ne a_1$ adjacent to both $b_5,b_i$. But then
$$\{\{a_2\},\{a_3,b_6\},\{a_4,b_1\},\{a_1,b_3\},\{b_4\},\{b_5,c_2,c_3,c_4,b_2\}\}$$
is a 6-cluster, a contradiction. This proves \ref{441}.~\bbox

\section{Excluding $K(3,4,0)$- and $K(2,6,0)$-subgraphs.}

Our main goal in this section is to eliminate $K(3,4,0)$-subgraphs; and to do this, we first eliminate
$K(3,7,3)$-subgraphs.

\begin{thm}\label{373}
Let $G$ be a minimal candidate with bipartition $(A,B)$. Then $G$ has no $K(3,7,3)$-subgraph.
\end{thm}
\Proof
Suppose that $a_1,a_2,a_3\in A$ are all adjacent to each of $b_1\LL b_7\in B$ except for the pairs $a_1b_1,a_2,b_2,a_3b_3$.
\\
\\
(1) {\em No vertex different from $a_1,a_2,a_3$ has
three neighbours in $\{b_4,b_5,b_6,b_7\}$.}
\\
\\
This is immediate since 
$G$ has no $K(4,4,1)$-subgraph by \ref{441}.
\begin{figure}[H]
\centering

\begin{tikzpicture}[scale=0.8,auto=left]
\tikzstyle{every node}=[inner sep=1.5pt, fill=black,circle,draw]
\node (a1) at (-3,2) {};
\node (a2) at (-1,2) {};
\node (a3) at (1,2) {};
\node (b1) at (1,0) {};
\node (b2) at (-1,0) {};
\node (b3) at (-3,0) {};
\node (b4) at (-4,4) {};
\node (b5) at (-2,4) {};
\node (b6) at (0,4) {};
\node (b7) at (2,4) {};

\foreach \from/\to in {a1/b2,a1/b3,a1/b4,a1/b5,a1/b6,a1/b7, a2/b1,a2/b3,a2/b4,a2/b5,a2/b6,a2/b7,a3/b1,a3/b2,a3/b4,
a3/b5, a3/b6,a3/b7}
\draw [-] (\from) -- (\to);
\tikzstyle{every node}=[]
\draw[left] (a1) node            {$a_1$};
\draw[left] (a2) node            {$a_2$};
\draw[right] (a3) node            {$a_3$};

\draw[below] (b1) node            {$b_1$};
\draw[below] (b2) node            {$b_2$};
\draw[below] (b3) node            {$b_3$};
\draw[above] (b4) node            {$b_4$};
\draw[above] (b5) node            {$b_5$};
\draw[above] (b6) node            {$b_6$};
\draw[above] (b7) node            {$b_7$};

\end{tikzpicture}
\caption{$K(3,7,3)$-subgraph.} \label{fig:373}
\end{figure}
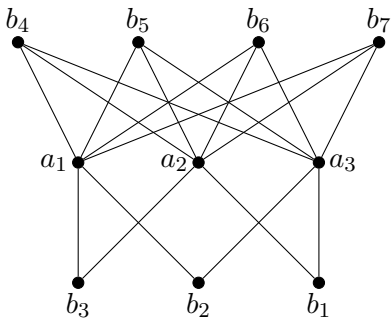
\noindent
(2) {\em No vertex different from $a_1,a_2,a_3$ has at least four neighbours in $\{b_1\LL b_7\}$.}
\\
\\
Suppose $a_4$ has at least four neighbours in $\{b_1\LL b_7\}$. Let $I$ be the set of $i\in \{1\LL 7\}$
such that $a_4,b_i$ are adjacent. Thus $|I|\ge 4$.
Let $b_8$ be a neighbour of $a_4$ not in $\{b_1\LL b_7\}$. Thus $b_8$ is nonadjacent to $a_1,a_2,a_3$, since
the latter have degree only six.
By \ref{mate}, for each $i\in I$ there exists 
$c_i\ne  a_1\LL a_4$ adjacent to $b_i, b_8$.

Suppose first that $1,2,3\in I$; and we may assume that $4\in I$ since $|I|\ge 4$.
Then 
$$\{\{a_4\},\{b_8,c_1,c_2,c_3,c_4\},\{a_1,b_2\},\{a_2,b_3\},\{a_3,b_1\},\{b_4\}\}$$
is a 6-cluster. So not all $1,2,3$ belong to $I$; and hence by (1), since $|I|\ge 4$, we may assume that $I=\{1,2,4,5\}$.
But then 
$$\{\{a_4\},\{b_8,c_1,c_2,c_4,c_5\},\{a_1,b_2\},\{a_2,b_3,b_4\},\{a_3,b_1\},\{b_5\}\}$$
is a 6-cluster. This proves (2).

\bigskip

Let $H$ be the cover graph with respect to $a_1,a_2,a_3,b_1\LL b_7$. By an argument like that in the proof of \ref{474}, it follows that
every partition of $V(H)$ into four sets is feasible,
and so $\chi(H)\ge 5$ by \ref{colouring}.  For each edge $b_ib_j$ of $H$ let $c_{i,j}\in A\setminus \{a_1,a_2,a_3\}$ be adjacent to
$b_i,b_j$.
\\
\\
(3) {\em The subgraph $H[\{b_4,b_5,b_6,b_7\}]$ has no triangle, and so is bipartite; and hence $b_1,b_2,b_3$ are pairwise adjacent in $H$.}
\\
\\
Suppose that say $b_4,b_5,b_6$ are pairwise adjacent in $H$. By (1), 
$c_{4,5}, c_{5,6}, c_{4,6}$ are all distinct. But then 
$$\{\{a_1,b_2\},\{a_2,b_3\},\{a_3,b_1\},\{c_{4,5},b_4\},\{c_{5,6},b_5\},\{c_{4,6},b_6\}\}$$
is a 6-cluster. 
So $H[\{b_4,b_5,b_6,b_7\}]$ is bipartite. Since
$\chi(H)\ge 5$ it follows that $b_1,b_2,b_3$ are pairwise adjacent in $H$. This proves (3).
\\
\\
(4) {\em Let $i\in \{1,2,3\}$, and let $j,k\in \{4,5,6,7\}$ be distinct. If $b_i, b_j, b_k$ are pairwise adjacent in $H$,
then $c_{i,j}, c_{i,k},c_{j,k}$ are all equal.}
\\
\\
Let $i = 3$, $j=4$ and $k=5$ say. 
Suppose first that $c_{3,4}, c_{4,5}, c_{3,5}$ are all different.
Since $c_{4,5}$ is different from $c_{2,3}$ by (2) it follows that
$$\{\{a_1\},\{a_2,b_6\},\{a_3,b_7\},\{b_3,c_{2,3},b_2,c_{3,4},c_{3,5}\},\{b_4\},\{c_{4,5},b_5\}\}$$
is a 6-cluster. Thus one of $c_{3,4}, c_{4,5}, c_{3,5}$ (say $c$) is adjacent to all of $b_3,b_4,b_5$.
Suppose that some $d\in \{c_{3,4}, c_{4,5}, c_{3,5}\}$ is different from $c$.
Then $d$ is different from one of $c_{1,3}, c_{2,3}$ by (2), say $d\ne c_{2,3}$,
and 
$$\{\{a_1\},\{a_2,b_6\},\{a_3,b_7\},\{b_3,c_{2,3},b_2,c\},\{b_4\},\{d,b_5\}\}$$
is a 6-cluster. 
This proves (4).
\\
\\
(5) {\em There exists a clique $X\subseteq V(H)$ of $H$ containing two of $b_1,b_2,b_3$ and two of $b_4\LL b_7$.}
\\
\\
If $H$ is perfect, then it has a clique of cardinality five, which therefore contains all of $b_1,b_2,b_3$ by (3)
and the claim holds.
Otherwise, $H$ has an odd hole or antihole as an induced subgraph; and since $H$ has only seven vertices and $\chi(H)\ge 5$,
it follows that $H$ has an induced cycle $C$ of length five, and the other two vertices of $H$ are adjacent to each other
and to every vertex of $C$. Since $b_1,b_2,b_3$ are pairwise adjacent, at least one of them is not in $V(C)$, say $b_1$;
and so at least three vertices of $C$ are not in $\{b_1,b_2,b_3\}$, and consequently an edge of $C$ has both ends
in $\{b_4\LL b_7\}$. But this set contains no triangle of $C$ by (3), and so the second vertex of $H$ not in $V(C)$ belongs
to $\{b_1,b_2,b_3\}$.
This proves (5).

\bigskip

From (5) we may assume that $b_2,b_3,b_4,b_5$ are pairwise adjacent in $H$.
By (4), $c_{1,4}, c_{1,5},c_{4,5}$ are all equal, and also 
$c_{1,5}, c_{1,6}, c_{5,6}$ are all equal. But then $c_{1,4}=c_{5,6}$ contrary to (2). 
This proves \ref{373}.~\bbox

\begin{thm}\label{340}
Let $G$ be a minimal candidate with bipartition $(A,B)$.
Then $G$ has no $K(3,4,0)$-subgraph.
\end{thm}
\Proof
Suppose that $a_1,a_2,a_3\in A$ are all adjacent to each of $b_1\LL b_4\in B$.
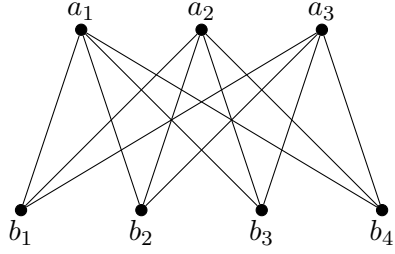
\begin{figure}[H]
\centering

\begin{tikzpicture}[scale=0.8,auto=left]
\tikzstyle{every node}=[inner sep=1.5pt, fill=black,circle,draw]
\node (a1) at (-2,3) {};
\node (a2) at (-0,3) {};
\node (a3) at (2,3) {};
\node (b1) at (-3,0) {};
\node (b2) at (-1,0) {};
\node (b3) at (1,0) {};
\node (b4) at (3,0) {};

\foreach \from/\to in {a1/b1,a1/b2,a1/b3,a1/b4, a2/b1,a2/b2,a2/b3,a2/b4,a3/b1,a3/b2,a3/b3,a3/b4}
\draw [-] (\from) -- (\to);
\tikzstyle{every node}=[]
\draw[above] (a1) node            {$a_1$};
\draw[above] (a2) node            {$a_2$};
\draw[above] (a3) node            {$a_3$};
\draw[below] (b1) node            {$b_1$};
\draw[below] (b2) node            {$b_2$};
\draw[below] (b3) node            {$b_3$};
\draw[below] (b4) node            {$b_4$};

\end{tikzpicture}
\caption{$K(3,4,0)$-subgraph.} \label{fig:340}
\end{figure}

No vertex in $A\setminus \{a_1,a_2,a_3\}$ has more than two neighbours in $\{b_1\LL b_4\}$, since $G$ has no $K(4,4,1)$-subgraph or 
$K(4,4,0)$-subgraph by \ref{441} and \ref{440}. Also
no vertex in $B\setminus \{b_1\LL b_4\}$ is adjacent to all three of $a_1,a_2,a_3$, since $G$ has no $K(3,5,0)$-subgraph by \ref{350}.

Let $B_0$ be the set of vertices in $B$ with exactly one neighbour in $\{a_1,a_2,a_3\}$. Thus $|B_0|$ is even.
\\
\\
(1) {\em $B_0\ne \emptyset$ and hence $|B_0|\ge 2$.}
\\
\\
Suppose that $B_0=\emptyset$. Since there are exactly six edges between $\{a_1,a_2,a_3\}$ and $B\setminus \{b_1\LL b_4\}$, it follows that there are exactly three vertices each adjacent to exactly two of $a_1,a_2,a_3$, and each 
of $a_1,a_2,a_3$ is adjacent to exactly two of these three vertices; but then $G$ contains a $K(3,7,3)$-subgraph, 
contrary to \ref{373}. This proves (1).
\\
\\
(2) {\em Every vertex adjacent to exactly two of $b_1\LL b_4$ is adjacent to every vertex in $B_0$.}
\\
\\
Let $a_4$ be adjacent to $b_1,b_2$ say, and let $b_5$ be adjacent to $a_1$ and not to $a_2,a_3$. Suppose
that $a_4,b_5$ are nonadjacent. By \ref{mate}, for $1\le i\le 4$ there is a vertex $c_i\ne a_1$ adjacent to 
$b_5,b_i$; and $c_i\ne a_2,a_3,a_4$ since these are not adjacent to $b_5$. But then 
$$\{\{a_4,b_1\},\{b_2\},\{b_5,c_1,c_2,c_3,c_4\},\{a_1\}, \{a_2,b_3\},\{a_3,b_4\}\}$$
is a 6-cluster. This proves (2).

\bigskip

Contracting $\{a_1,a_2,a_3,b_1,b_2,b_3,b_4\}$ into $B$ does not yield a smaller candidate, so some vertex $a_4$ 
different from $a_1,a_2,a_3$ has at least two (and hence exactly two) neighbours in $\{b_1,b_2,b_3,b_4\}$. From
the symmetry we may assume that $a_4$ is adjacent to $b_1,b_2$. Since $a_4$ has degree six, and is adjacent to every
vertex in $B_0$, it follows that $|B_0|\le 4$, and so some vertex is adjacent to exactly two of $a_1,a_2,a_3$.
\\
\\
(3) {\em There is a vertex different from $a_1\LL a_4$ adjacent to $b_3,b_4$.}
\\
\\
Suppose not. Choose distinct $b_5,b_6\in B_0$. 
By \ref{mate}, for $i\in \{3,4\}$ and $j\in \{5,6\}$ there is a vertex $c_{i,j}\notin \{a_1\LL a_4\}$ 
adjacent to
$b_i,b_j$; and $c_{i,j}\ne a_4$ since $a_4$
has only two neighbours in $\{b_1\LL b_4\}$. Moreover $\{c_{3,5},c_{3,6}\}$ is disjoint from $\{c_{4,5}, c_{4,6}\}$,
since these vertices only have one neighbour in $\{b_3,b_4\}$ by hypothesis. Let $b_7$ be a vertex adjacent to exactly two of $a_1,a_2,a_3$, say $a_2,a_3$; then since $a_4$ is adjacent to $b_5,b_6$ by (2),
$$\{\{a_1,b_1\},\{a_2\},\{a_3,b_7\},\{c_{3,5},c_{3,6},b_3,b_5\},\{c_{4,5},c_{4,6},b_4\},\{a_4,b_2,b_6\}\}$$
is a 6-cluster. This proves (3).

\begin{figure}[H]
\centering

\begin{tikzpicture}[scale=0.8,auto=left]
\tikzstyle{every node}=[inner sep=1.5pt, fill=black,circle,draw]
\node (a1) at (-2,3) {};
\node (a2) at (-0,3) {};
\node (a3) at (2,3) {};
\node (a4) at (-2,-2) {};
\node (a5) at (2,-2) {};
\node (b1) at (-3,0) {};
\node (b2) at (-1,0) {};
\node (b3) at (1,0) {};
\node (b4) at (3,0) {};

\foreach \from/\to in {a1/b1,a1/b2,a1/b3,a1/b4, a2/b1,a2/b2,a2/b3,a2/b4,a3/b1,a3/b2,a3/b3,a3/b4,a4/b1,a4/b2,a5/b3,a5/b4}
\draw [-] (\from) -- (\to);
\tikzstyle{every node}=[]
\draw[above] (a1) node            {$a_1$};
\draw[above] (a2) node            {$a_2$};
\draw[above] (a3) node            {$a_3$};
\draw[below] (a4) node            {$a_4$};
\draw[below] (a5) node            {$a_5$};
\draw[left] (b1) node            {$b_1$};
\draw[right] (b2) node            {$b_2$};
\draw[left] (b3) node            {$b_3$};
\draw[right] (b4) node            {$b_4$};

\end{tikzpicture}
\caption{For the proof of \ref{340}.} \label{fig:340+}
\end{figure}
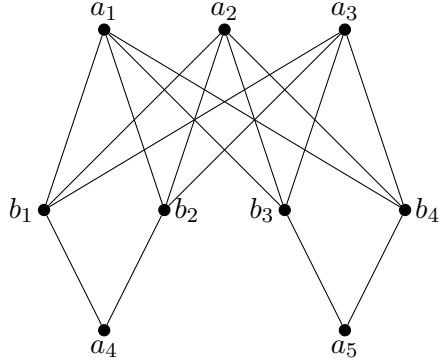

\bigskip

Let $a_5\ne a_1\LL a_4$ be adjacent to $b_3,b_4$. By (2) $a_5$ is adjacent to every vertex in $B_0$.
\\
\\
(4) {\em There is no path $P$ of $G\setminus \{a_1\LL a_5\}$ with ends in distinct sets in the list $\{b_1\},\{b_2\},\{b_3\},\{b_4\}, B_0$.}
\\
\\
Suppose that $P$ is such a path. 
By choosing $P$ minimal we may assume that no internal vertex of $P$ belongs to
$\{b_1,b_2,b_3,b_4\}\cup B_0$; and from the symmetry we may assume that $b_1$ is an end of $P$. Let the other end be $b_i$ say.
Up to symmetry there are three cases: $i = 2$, $i=3$,  and $i=5$ for some $b_5\in B_0$, and in the third case we may assume that 
$a_1$ is adjacent to $b_5$ from the symmetry. 
\begin{itemize}
\item If $i=2$ then $\{\{a_1\},\{a_2,b_3\},\{a_3,b_4\},\{b_1\}, P^*\cup \{b_2\},B_0\cup \{a_4,a_5\}\}$ is a 6-cluster.
\item If $i=3$ then $\{\{a_1\},\{a_2,b_2\},\{a_3,b_4\},\{b_1\}, P^*\cup \{b_3\},B_0\cup \{a_4,a_5\}\}$ is a 6-cluster.
\item If $i=5$ then $\{\{a_1\},\{a_2, b_3\},\{a_3,b_4\},P^*\cup \{b_1\},\{a_4,b_2\},\{a_5,b_6\}\}$ is a 6-cluster.
\end{itemize}
This proves (4).
\\
\\
(5) {\em $b_1\LL b_4$ all have degree four in $G$.}
\\
\\
From (4), $b_1,b_2,b_3,b_4,b_5$ all belong to different components of $G\setminus \{a_1\LL a_5\}$, and none of these four
compoents contains any vertex of $B_0$. By 
\ref{facts}, $G\setminus \{a_1\LL a_5\}$ has exactly five components, and four of them are singletons;
and therefore one contains all of $B_0$ and so is not a singleton.
This proves (5).

\bigskip

As we observed earlier, there is a vertex $c\in B$ that is adjacent to exactly two of $a_1,a_2,a_3$, say to $a_2,a_3$.
Moreover, since $B_0\ne \emptyset$, not both neighbours of $a_1$ in $B\setminus \{b_1\LL b_4\}$
have a neighbour in $\{b_2,b_3\}$; and so $a_1$ has a neighbour in $B_0$, say $b_5$. Since we cannot obtain a smaller candidate
by contracting $\{a_1\LL a_5, b_1\LL b_5,c\}$ into $B$, there exists $a_6\in A$ different from $a_1\LL a_5$ with two neighbours in
$\{b_1\LL b_5,c\}$. By (5), $a_6$ is adjacent to $b_5,c$. Let $b_6\in B_0\setminus \{b_5\}$. Then 
$$\{\{a_1\},\{a_2,b_1\},\{a_3,b_3\},\{a_4,b_2\},\{a_5,b_4\},\{a_6,c,b_5\}\}$$
is a 6-cluster. This proves \ref{340}.~\bbox

\section{Excluding $K(2,5,0)$-subgraphs}

Our next goal is to eliminate $K(2,5,0)$-subgraphs. We begin with:
\begin{thm}\label{260}
Let $G$ be a minimal candidate with bipartition $(A,B)$. Then $G$ has no $K(2,6,0)$-subgraph.
\end{thm}
\Proof
Suppose that $a_1,a_2\in A$ are both adjacent to each of $b_1\LL b_6\in B$. The cover graph $H$ with respect to
$a_1,a_2, b_1\LL b_6$ has chromatic number at least five, by \ref{colouring} (note that every partition of $\{b_1\LL b_6\}$
into four sets is feasible, since at least two of them will be singletons and therefore already induce connected subgraphs).
Consequently $H$ has a clique of size five, and so we may assume that $b_1\LL b_5$ are pairwise adjacent in $H$.
By \ref{340}, no vertex in $A\setminus \{a_1,a_2\}$ has more than three neighbours in $\{b_1\LL b_6\}$.
\begin{figure}[H]
\centering

\begin{tikzpicture}[scale=0.8,auto=left]
\tikzstyle{every node}=[inner sep=1.5pt, fill=black,circle,draw]
\node (a1) at (-2,3) {};
\node (a2) at (2,3) {};
\node (b1) at (-5,0) {};
\node (b2) at (-3,0) {};
\node (b3) at (-1,0) {};
\node (b4) at (1,0) {};
\node (b5) at (3,0) {};
\node (b6) at (5,0) {};

\foreach \from/\to in {a1/b1,a1/b2,a1/b3,a1/b4, a1/b5, a1/b6,a2/b1,a2/b2,a2/b3,a2/b4,a2/b5, a2/b6}
\draw [-] (\from) -- (\to);
\tikzstyle{every node}=[]
\draw[above] (a1) node            {$a_1$};
\draw[above] (a2) node            {$a_2$};

\draw[below] (b1) node            {$b_1$};
\draw[below] (b2) node            {$b_2$};
\draw[below] (b3) node            {$b_3$};
\draw[below] (b4) node            {$b_4$};
\draw[below] (b5) node            {$b_5$};
\draw[below] (b6) node            {$b_6$};

\end{tikzpicture}
\caption{$K(2,6,0)$-subgraph.} \label{fig:260}
\end{figure}
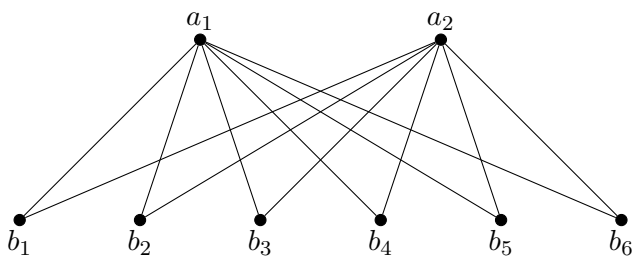

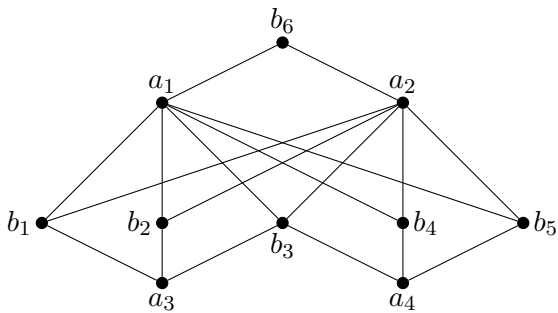
\begin{figure}[H]
\centering

\begin{tikzpicture}[scale=0.8,auto=left]
\tikzstyle{every node}=[inner sep=1.5pt, fill=black,circle,draw]
\node (a1) at (-2,2) {};
\node (a2) at (2,2) {};
\node (b1) at (-4,0) {};
\node (b2) at (-2,0) {};
\node (b3) at (-0,0) {};
\node (b4) at (2,0) {};
\node (b5) at (4,0) {};
\node (b6) at (0,3) {};
\node (a3) at (-2,-1) {};
\node (a4) at (2,-1) {};

\foreach \from/\to in {a1/b1,a1/b2,a1/b3,a1/b4, a1/b5, a1/b6,a2/b1,a2/b2,a2/b3,a2/b4,a2/b5, a2/b6, a3/b1,a3/b2,a3/b3,a4/b3,a4/b4,a4/b5}
\draw [-] (\from) -- (\to);
\tikzstyle{every node}=[]
\draw[above] (a1) node            {$a_1$};
\draw[above] (a2) node            {$a_2$};

\draw[left] (b1) node            {$b_1$};
\draw[left] (b2) node            {$b_2$};
\draw[below] (b3) node            {$b_3$};
\draw[right] (b4) node            {$b_4$};
\draw[right] (b5) node            {$b_5$};
\draw[above] (b6) node            {$b_6$};

\draw[below] (a3) node            {$a_3$};
\draw[below] (a4) node            {$a_4$};

\end{tikzpicture}
\caption{For the proof of \ref{260}.} \label{fig:260+}
\end{figure}

For $1\le i<j\le 5$ let $c_{i,j}\in A\setminus \{a_1,a_2\}$ be adjacent to $b_i, b_j$. If the six vertices
$c_{i,j}\;(1\le i<j\le 4)$
are all distinct, there is a 6-cluster
$$\{\{b_1,c_{1,2}, c_{1,3}, c_{1,4}\},\{b_2,c_{2,3}, c_{2,4}\},\{b_3,c_{3,4}\},\{b_4\},\{a_1\},\{a_2,b_5\}\},$$
a contradiction. So we may assume that some two are equal, and hence some vertex in $A\setminus \{a_1,a_2\}$ is adjacent
to three of $b_1,b_2,b_3,b_4$; say $a_3$ is adjacent to $b_1,b_2,b_3$.
If none of the vertices $c_{i,j}\;(i\in \{1,2,3\}, j\in \{4,5\})$
is adjacent to both $b_4,b_5$, then
$$\{\{b_1\},\{b_2, a_3\},\{b_4,c_{1,4}, c_{2,4}\},\{b_5,c_{1,5},c_{2,5},c_{4,5}\},\{a_1\},\{a_2,b_6\}\},$$
is a 6-cluster, a  contradiction; so we may assume that some $a_4$ is adjacent to $b_3,b_4,b_5$ say. (See figure \ref{fig:260+}.)

If some $a_5$ different from $a_1\LL a_4$ is adjacent to three of $b_1,b_2,b_4,b_5$, say to $b_1,b_2, b_4$, then
$$\{\{b_1\},\{b_2, a_3\},\{b_4,a_5\},\{b_5, a_4,c_{1,5}, c_{2,5}\}, \{a_1\}, \{a_2,b_6\}\}$$
is a 6-cluster, a  contradiction; so we may assume that no vertex different from $a_1\LL a_4$ is adjacent to three of
$b_1,b_2, b_4, b_5$. Consequently $ c_{1,4}, c_{2,4},  c_{1,5}, c_{2,5}$ are all different. But then
$$\{\{b_1\},\{b_2, a_3\},\{b_4, c_{1,4}, c_{2,4}\},\{b_5, a_4, c_{1,5}, c_{2,5}\}, \{a_1\}, \{a_2,b_6\}\}$$
is a 6-cluster, a contradiction. This proves \ref{260}.~\bbox

\begin{thm}\label{250}
Let $G$ be a minimal candidate with bipartition $(A,B)$. Then $G$ has no $K(2,5,0)$-subgraph.
\end{thm}
\Proof
Suppose that $a_1,a_2\in A$ are both adjacent to each of $b_1\LL b_5\in B$. 
Let $b_6, b_7$ be the neighbours of $a_1,a_2$ respectively
that are not in $\{b_1\LL b_5\}$. Thus $b_6\ne b_7$ since there is no $K(2,6,0)$-subgraph by \ref{260}. No vertex different from
$a_1,a_2$ has four neighbours in $\{b_1\LL b_5\}$ since there is no $K(3,4,0)$-subgraph by \ref{340}.
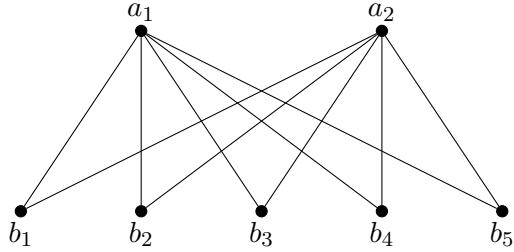
\begin{figure}[H]
\centering

\begin{tikzpicture}[scale=0.8,auto=left]
\tikzstyle{every node}=[inner sep=1.5pt, fill=black,circle,draw]
\node (a1) at (-2,3) {};
\node (a2) at (2,3) {};
\node (b1) at (-4,0) {};
\node (b2) at (-2,0) {};
\node (b3) at (-0,0) {};
\node (b4) at (2,0) {};
\node (b5) at (4,0) {};

\foreach \from/\to in {a1/b1,a1/b2,a1/b3,a1/b4, a1/b5, a2/b1,a2/b2,a2/b3,a2/b4,a2/b5}
\draw [-] (\from) -- (\to);
\tikzstyle{every node}=[]
\draw[above] (a1) node            {$a_1$};
\draw[above] (a2) node            {$a_2$};

\draw[below] (b1) node            {$b_1$};
\draw[below] (b2) node            {$b_2$};
\draw[below] (b3) node            {$b_3$};
\draw[below] (b4) node            {$b_4$};
\draw[below] (b5) node            {$b_5$};

\end{tikzpicture}
\caption{$K(2,5,0)$-subgraph.} \label{fig:250}
\end{figure}

The cover graph $H$ with respect to
$a_1,a_2, b_1\LL b_5$ has chromatic number at least four, by \ref{colouring},
and so has a clique of cardinality four, say $b_1,b_2,b_3,b_4$. 
For all distinct $i,j\in \{1,2,3,4,6,7\}$ 
let $c_{i,j}\in A\setminus \{a_1,a_2\}$ be adjacent to $b_i, b_j$, if there is such a vertex. Thus $c_{i,j}$ exists 
for $1\le i<j\le 4$, and also for all $i\in \{1\LL 5\}$ and $j\in \{6,7\}$, by \ref{mate}.
\\
\\
(1) {\em Some vertex in $A\setminus \{a_1,a_2\}$ has three neighbours in $\{b_1\LL b_4\}$.}
\\
\\
Suppose not. Then each of the vertices $c_{i,j}\; (1\le i<j\le 4)$ has only two neighbours in $\{b_1\LL b_4\}$,
and in particular they are all different. But then
$$\{\{b_1\},\{c_{1,2},b_2\},\{c_{1,3},c_{2,3},b_3\},\{c_{1,4},c_{2,4},c_{3,4},b_4\},\{a_1\},\{a_2,b_5\}\}$$
is a 6-cluster. This proves (1).

\bigskip

Thus we may assume that some $a_3$ is adjacent to $b_1,b_2,b_3$.
\\
\\
(2) {\em Some vertex in  $A\setminus \{a_1,a_2, a_3\}$ is adjacent to $b_4$ and to two of $b_1,b_2,b_3$.}
\\
\\
Suppose not; so $c_{1,4}, c_{2,4}, c_{3,4}$ are all different. Suppose that some vertex $c\ne a_1,a_2,a_3$ is adjacent to two 
of $b_1,b_2,b_3$,
say to $b_1,b_2$. Then 
$$\{\{b_1\},\{c,b_2\},\{a_3,b_3\},\{c_{1,4},c_{2,4},c_{3,4},b_4\},\{a_1\},\{a_2,b_5\}\}$$
is a 6-cluster. Thus there is no such $c$. There is a vertex $b_8\in B\setminus \{b_1\LL b_7\}$ adjacent to 
$a_3$; and by \ref{mate}, for $i = 1,2,3$ there exists $d_i\in A\setminus \{a_1,a_2,a_3\}$ adjacent to $b_8,b_i$. 
Since no $c\ne a_1,a_2,a_3$ is adjacent to two 
of $b_1,b_2,b_3$, it follows that $d_1,d_2,d_3$ are all different. Consequently $c_{i,4}\ne d_j$ for all distinct $i,j\in \{1,2,3\}$.
But then 
$$\{\{c_{1,4},d_1,b_1\},\{a_3,c_{2,4},d_2,b_2\},\{c_{3,4},d_3, b_3,b_8\}, \{a_1\},\{a_2,b_5\}\}$$
is a 6-cluster. This proves (2).

\bigskip
\begin{figure}[H]
\centering

\begin{tikzpicture}[scale=0.8,auto=left]
\tikzstyle{every node}=[inner sep=1.5pt, fill=black,circle,draw]
\node (a1) at (-2,2) {};
\node (a2) at (2,2) {};
\node (b1) at (-3,0) {};
\node (b2) at (-1,0) {};
\node (b3) at (1,0) {};
\node (b4) at (3,0) {};
\node (b5) at (0,3) {};
\node (a3) at (-2,-2) {};
\node (a4) at (2,-2) {};
\node (a5) at (0,-2) {};
\node (b6) at (-5,0) {};
\node (b7) at (5,0) {};

\foreach \from/\to in {a1/b1,a1/b2,a1/b3,a1/b4, a1/b5, a2/b1,a2/b2,a2/b3,a2/b4,a2/b5,a3/b1,a3/b2,a3/b3,a4/b2,a4/b3,a4/b4, a5/b1,a5/b4,a1/b6,a2/b7}
\draw [-] (\from) -- (\to);
\tikzstyle{every node}=[]
\draw[above] (a1) node            {$a_1$};
\draw[above] (a2) node            {$a_2$};

\draw[left] (b1) node            {$b_1$};
\draw[left] (b2) node            {$b_2$};
\draw[right] (b3) node            {$b_3$};
\draw[right] (b4) node            {$b_4$};
\draw[above] (b5) node            {$b_5$};
\draw[below] (a3) node            {$a_3$};
\draw[below] (a4) node            {$a_4$};
\draw[below] (a5) node            {$a_5$};
\draw[left] (b6) node            {$b_6$};
\draw[right] (b7) node            {$b_7$};

\end{tikzpicture}
\caption{For the proof of \ref{250}.} \label{fig:250+}
\end{figure}
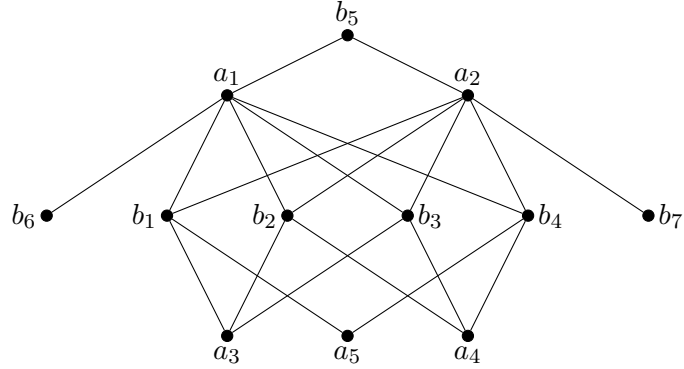

Thus we may assume that $a_4\in \setminus \{a_1,a_2,a_3\}$ is adjacent to $b_2,b_3,b_4$.
Since $\{b_1,b_2,b_3,b_4\}$ is a clique of $H$, there exists $a_5\in  A\setminus \{a_1,a_2,a_3,a_4\}$ adjacent to $b_1,b_4$.
(See figure \ref{fig:250+}.)
\\
\\
(3) {\em No vertex in $A$ different from $a_1,a_2,a_3,a_4$ has two neighbours in $\{b_1,b_2,b_3\}$, or has two neighbours in $\{b_2,b_3,b_4\}$.}
\\
\\
Suppose that $a_6$ is such a vertex, adjacent to two of $b_1,b_2,b_3$ say. (Possibly $a_6=a_5$.)
If $a_6$ is adjacent to $b_2,b_3$, then $a_6\ne a_5$, and
$$\{\{a_3,b_1\},\{b_2\},\{a_6,b_3\},\{a_4,a_5,b_4\},\{a_1\},\{a_2,b_5\}\}$$
is a 6-cluster. Thus from the symmetry we may assume that $a_6$ is adjacent to $b_1,b_2$, and now possibly $a_6=a_5$.
Then 
$$\{\{a_5,a_6,b_1\},\{b_2\},\{a_3,b_3\},\{a_4,b_4\},\{a_1\},\{a_2,b_5\}\}$$
is a 6-cluster. This proves (3).
\\
\\
(4) {\em Every vertex in $B\setminus \{b_1\LL b_7\}$ adjacent to one of $a_3,a_4$ is adjacent to both $a_3,a_4$.}
\\
\\
Suppose that $b_8\in B\setminus \{b_1\LL b_7\}$ is adjacent to $b_3$ and not to $b_4$.
By \ref{mate}, for $i=2,3$
there exists $c_{i,8}\in A\setminus \{a_3\}$ adjacent to $b_i, b_8$; $c_{i,8}\ne a_1,a_2,a_4$ since $a_1,a_2,a_4$
are not adjacent to $b_8$, and $c_{i,8}\ne a_5$ since $a_5$ is not adjacent to $b_i$ by (3).
But then
$$\{\{a_3,b_1\},\{c_{2,8},c_{3,8},b_8,b_2\}, \{b_3\},\{a_4,a_5,b_4\},\{a_1\},\{a_2,b_5\}\}$$
is a 6-cluster. This proves (4).
\\
\\
(5) {\em Each of $b_6,b_7$ is adjacent to at least one of $a_3,a_4$.}
\\
\\
Suppose that $b_6$ is nonadjacent to both $a_3,a_4$. Thus $c_{i,6}\ne a_4,a_5$ for $1\le i\le 4$. But then 
$$\{\{c_{1,6},c_{4,6},b_1,b_6\},\{a_3,b_2\},\{a_4,b_3\},\{b_4\},\{a_1\},\{a_2,b_5\}\}$$
is a 6-cluster. This proves (5).
\\
\\
(6) {\em Each of $b_6,b_7$ is adjacent to both of $a_3,a_4$.}
\\
\\
\begin{figure}[H]
\centering

\begin{tikzpicture}[scale=0.8,auto=left]
\tikzstyle{every node}=[inner sep=1.5pt, fill=black,circle,draw]
\node (a1) at (-2,2) {};
\node (a2) at (2,2) {};
\node (b1) at (-3,0) {};
\node (b2) at (-1,0) {};
\node (b3) at (1,0) {};
\node (b4) at (3,0) {};
\node (b5) at (0,3) {};
\node (a3) at (-2,-2) {};
\node (a4) at (2,-2) {};
\node (a5) at (0,-2) {};
\node (b6) at (-5,0) {};
\node (b7) at (5,0) {};
\node (b8) at (0,-3) {};

\foreach \from/\to in {a1/b1,a1/b2,a1/b3,a1/b4, a1/b5, a2/b1,a2/b2,a2/b3,a2/b4,a2/b5,a3/b1,a3/b2,a3/b3,a3/b6,a4/b2,a4/b3,a4/b4, a4/b7,a5/b1,a5/b4,a1/b6,a2/b7,a3/b8,a4/b8}
\draw [-] (\from) -- (\to);
\tikzstyle{every node}=[]
\draw[above] (a1) node            {$a_1$};
\draw[above] (a2) node            {$a_2$};

\draw[left] (b1) node            {$b_1$};
\draw[left] (b2) node            {$b_2$};
\draw[right] (b3) node            {$b_3$};
\draw[right] (b4) node            {$b_4$};
\draw[above] (b5) node            {$b_5$};
\draw[below] (a3) node            {$a_3$};
\draw[below] (a4) node            {$a_4$};
\draw[below] (a5) node            {$a_5$};
\draw[left] (b6) node            {$b_6$};
\draw[right] (b7) node            {$b_7$};
\draw[below] (b8) node            {$b_8$};

\end{tikzpicture}
\caption{For the proof of \ref{250}, step (6).} \label{fig:250+6}
\end{figure}

Suppose that $a_4,b_6$ are nonadjacent. By (5), $a_3b_6$ is an edge.
From (4), since $a_3,a_4$ both have degree six, and have the same number of neighbours in $B\setminus \{b_1\LL b_7\}$,
it follows that they have the same number of neighbours in $\{b_6,b_7\}$; and so $a_4$ has at least one neighbour in $\{b_6,b_7\}$,
and therefore $a_4b_7$ is an edge, and $a_3, b_7$ are not adjacent. Suppose that $c_{4,6}, c_{1,7}$ are different; then from the symmetry we may assume that
$c_{1,7}\ne a_5$.
Then
$$\{\{c_{1,7},b_1\},\{a_3,b_2,b_6\},\{a_4, b_3,b_7\},\{a_5,c_{4,6},b_4\},\{a_1\},\{a_2,b_5\}\}$$
is a 6-cluster. So $c_{1,7}=c_{4,6}$, and we may assume that they both equal $a_5$. Since $a_5$ has at most five neighbours in
$\{b_1\LL b_7\}$, it has a neighbour $b_8$ different from $b_1\LL b_7$. Consequently $b_8$ is nonadjacent to $a_1,a_2$.
Suppose that $b_8$ is nonadjacent to both $a_3,a_4$.
By \ref{mate}, for $i = 1,4,6,7$ there exists $c_{i,8}$ adjacent to $b_i, b_8$, and different from $a_5$. Since $b_8$
is nonadjacent to $a_1,a_2,a_3,a_4$, it follows that $c_{i,8}\ne a_1\LL a_5$. But then
$$\{\{a_5,b_1\},\{a_3,b_2\},\{a_4,b_3,b_7\},\{c_{4,8},c_{6,8},b_4,b_6,b_8\},\{a_1\},\{a_2,b_5\}\}$$
is a 6-cluster. So $b_8$ is adjacent to one of $a_3,a_4$, and hence to both $a_3,a_4$ by (2). By \ref{mate},
there is a vertex $a_6\ne a_5$ adjacent to $b_6,b_7$; and so $a_6\ne a_1,a_2,a_3,a_4$ since none of these four vertices is 
adjacent to both $b_6,b_7$. 
But then
$$\{\{a_1,b_1\},\{a_3,b_2\},\{a_4,b_3\},\{a_5,b_4\},\{a_2\},\{a_6,b_6,b_7\}\}$$
is a 6-cluster. This proves (6).

\bigskip

\begin{figure}[H]
\centering

\begin{tikzpicture}[scale=0.8,auto=left]
\tikzstyle{every node}=[inner sep=1.5pt, fill=black,circle,draw]
\node (a1) at (-3,0) {};
\node (a2) at (-1,0) {};
\node (b1) at (-3,-2) {};
\node (b2) at (-1,2) {};
\node (b3) at (1,2) {};
\node (b4) at (-1,-2) {};
\node (b5) at (-3,2) {};
\node (a3) at (1,0) {};
\node (a4) at (3,0) {};
\node (b6) at (1,-2) {};
\node (b7) at (3,-2) {};
\node (b8) at (3,2) {};

\foreach \from/\to in {a1/b1,a1/b2,a1/b3,a1/b4, a1/b5, a2/b1,a2/b2,a2/b3,a2/b4,a2/b5,a3/b1,a3/b2,a3/b3,a3/b6,a3/b7,a4/b2,a4/b3,a4/b4,a4/b6,a4/b7,a1/b6,a2/b7,a3/b8,a4/b8}
\draw [-] (\from) -- (\to);
\tikzstyle{every node}=[]
\draw[left] (a1) node            {$a_1$};
\draw[left] (a2) node            {$a_2$};

\draw[below] (b1) node            {$b_1$};
\draw[above] (b2) node            {$b_2$};
\draw[above] (b3) node            {$b_3$};
\draw[below] (b4) node            {$b_4$};
\draw[above] (b5) node            {$b_5$};
\draw[right] (a3) node            {$a_3$};
\draw[right] (a4) node            {$a_4$};
\draw[below] (b6) node            {$b_6$};
\draw[below] (b7) node            {$b_7$};
\draw[above] (b8) node            {$b_8$};

\end{tikzpicture}
\caption{For the last part of the proof of \ref{250}. ($a_5$ is not drawn.)} \label{fig:250last}
\end{figure}
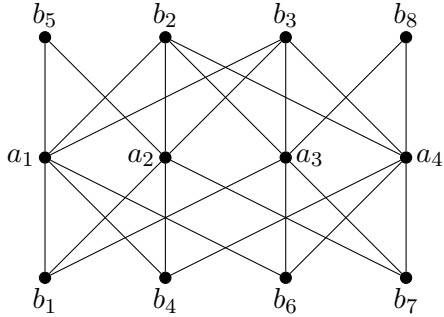

From (6), there is a unique vertex $b_8\in B\setminus \{b_1\LL b_7\}$ adjacent to $a_3$, and it is adjacent to 
both of $a_3,a_4$ by (4).
The subgraph induced on $\{a_1,a_2,a_3,a_4,b_1,b_2,b_3,b_4,b_5,b_6,b_7,b_8\}$ (note that $a_5$ is not included) has some significant 
symmetry, which will help reduce the case analysis to come. There are symmetries that exchange 
\begin{itemize}
\item $b_2$ with $b_3$;
\item $a_1$ with $a_2$, and $b_6$ with $b_7$; 
\item $a_3$ with $a_4$, and $b_1$ with $b_4$;
\item $a_1$ with $a_3$, $a_2$ with $a_4$, $b_5$ with $b_8$, $b_6$ with $b_1$, and $b_7$ with $b_4$.
\end{itemize}
Let us call these symmetries the first, second, third and fourth symmetries respectively. In the argument to come, we will avoid 
making use of $a_5$, in order to maintain these symmetries. 
\\
\\
(7) {\em Let $C_2$ be the component of $G\setminus \{a_1,a_2,a_3,a_4\}$ that contains $b_2$; then it contains none of 
$b_1\LL b_8$ except $b_2$.}
\\
\\
Suppose not. 
By \ref{mate}, $c_{5,6},c_{5,7},c_{1,8},c_{4,8}$ exist, and they are
different from $a_1\LL a_4$. Let $X=\{c_{5,6},c_{5,7}, b_5,b_6,b_7\}$ and $Y=\{c_{1,8}, c_{4,8}, b_1,b_4,b_8\}$. ($X,Y$ might not 
be disjoint.) Since $C_1$ contains one of $b_1,b_3\LL b_8$, there is a 
a minimal path $P$ of $G\setminus \{a_1,a_2,a_3,a_4\}$ with one end $b_2$ and the other end in $\{b_3\}\cup X\cup Y$.
From the minimality of $P$ it has no other neighbour in $\{b_3\}\cup X\cup Y$. Let $z$ be the end of $P$ different from $b_2$.
If $z=b_3$ then 
$$\{\{b_3\},P^*\cup \{b_2\},\{a_1\},\{a_2,b_5\},\{a_3,b_1\},\{a_4,b_4\}\}$$ 
is a 6-cluster. Thus $z\in X\cup Y$, and from the fourth symmetry we may assume that $z\in X$, and from the second symmetry we
may assume that $z\in \{a_1,b_5,b_6\}$.
But then 
$$\{P^*\cup X, \{a_1,b_4\},\{a_2\}, \{a_3,b_1\}, \{b_2\}, \{a_4,b_3\}\}$$
is a 6-cluster. This proves (7).

\bigskip

Similarly, let $C_3$ be the component of $G\setminus \{a_1,a_2,a_3,a_4\}$ that contains $b_3$; then it contains none of
$b_1\LL b_8$ except $b_3$. Now $a_1\LL a_4$ are the only vertices in $A\setminus V(C_i)$ that have a neighbour in $V(C_i)$,
for $i = 2,3$, and in particular, no vertex in $V(C_i)\cap A$ has a neighbour in $V(G)\setminus V(C_i)$. Since
$C_i$ is not a smaller candidate, it follows that $|A\cap V(C_i)|<|B\cap V(C_i)|$, for $i = 2,3$. But there are at least six
vertices in $B\setminus (V(C_2)\cup V(C_3))$ that have a neighbour in $\{a_1,a_3\}$, namely $b_1,b_4,b_5,b_6,b_7,b_8$;
and similarly there are six such vertices that have a neighbour in $\{a_2,a_4\}$. Hence contracting $V(C_2)\cup \{a_1,a_3\}$
and $V(C_2)\cup \{a_2,a_4\}$ into $A$ makes a smaller candidate, a contradiction. This proves \ref{250}.~\bbox

\section{Excluding $K(4,4,2)$-subgraphs}

\begin{thm}\label{442}
Let $G$ be a minimal candidate with bipartition $(A,B)$. Then $G$ has no $K(4,4,2)$-subgraph.
\end{thm}
\Proof
Suppose that $a_1\LL a_4\in A$ are all adjacent to each of $b_1\LL b_4\in B$ except the pairs $a_1b_1,a_2b_2$.
\begin{figure}[H]
\centering

\begin{tikzpicture}[scale=0.8,auto=left]
\tikzstyle{every node}=[inner sep=1.5pt, fill=black,circle,draw]
\node (a1) at (-3,2) {};
\node (a2) at (3,2) {};
\node (a3) at (-1,2) {};
\node (a4) at (1,2) {};
\node (b1) at (-3,-1) {};
\node (b2) at (3,-1) {};
\node (b3) at (-1,-1) {};
\node (b4) at (1,-1) {};

\foreach \from/\to in {a1/b2,a1/b3,a1/b4, a2/b1,a2/b3,a2/b4,a3/b1,a3/b2,a3/b3,a3/b4,a4/b1,a4/b2,a4/b3,a4/b4}
\draw [-] (\from) -- (\to);
\tikzstyle{every node}=[]
\draw[above] (a1) node            {$a_1$};
\draw[above] (a2) node            {$a_2$};
\draw[above] (a3) node            {$a_3$};
\draw[above] (a4) node            {$a_4$};
\draw[below] (b1) node            {$b_1$};
\draw[below] (b2) node            {$b_2$};
\draw[below] (b3) node            {$b_3$};
\draw[below] (b4) node            {$b_4$};

\end{tikzpicture}
\caption{$K(4,4,2)$-subgraph.} \label{fig:442}
\end{figure}
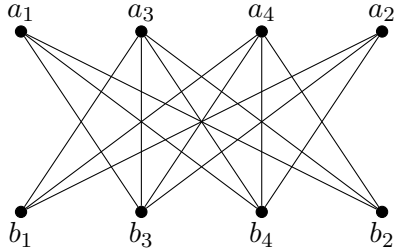

For $i = 3,4$, $a_i$ has two neighbours, $d_{1,i}, d_{2,i}$ say, not in $\{b_1\LL b_4\}$. 
Since there is no $K(2,5,0)$-subgraph by \ref{250},
none of $d_{1,3}, d_{2,3}, d_{1,4}, d_{2,4}$ is adjacent to both $a_3,a_4$, and so they are all distinct.
\\
\\
(1) {\em Each of $d_{1,3}, d_{2,3}, d_{1,4}, d_{2,4}$ has a neighbour in $\{a_1,a_2\}$.}
\\
\\
Suppose that $d_{1,3}$ say is nonadjacent to both $a_1,a_2$. By \ref{mate}, for $1\le i\le 4$ there is a vertex $c_i$ different from
$a_3$ that is adjacent to both $d_{1,3},b_i$. But then 
$$\{\{a_3\},\{c_1,c_2,c_3,c_4,d_{1,3}\},\{a_1,b_2\},\{a_2,b_3\},\{a_3,b_1\},\{b_4\}\}$$
is a 6-cluster. This proves (1).

\bigskip

Not both $d_{1,3}, d_{2,3}$ are adjacent to $a_1$, since there is no $K(2,5,0)$-subgraph by \ref{250}, and similarly they are not 
both adjacent to $a_2$, so we may assume that $d_{1,3}$
is nonadjacent to $a_2$, and $d_{2,3}$ is nonadjacent to $a_1$. The same applies for $d_{1,4}, d_{2,4}$; so for each $i\in \{1,2\}$
and each $j\in \{3,4\}$, $d_{i,j}$ is adjacent to $a_i, a_j$ and nonadjacent to the other two vertices in $\{a_1\LL a_4\}$.
\begin{figure}[H]
\centering

\begin{tikzpicture}[scale=0.8,auto=left]
\tikzstyle{every node}=[inner sep=1.5pt, fill=black,circle,draw]
\node (a1) at (-3,2) {};
\node (a2) at (3,2) {};
\node (a3) at (-1,2) {};
\node (a4) at (1,2) {};
\node (b1) at (-3,-1) {};
\node (b2) at (3,-1) {};
\node (b3) at (-1,-1) {};
\node (b4) at (1,-1) {};
\node (d13) at (-3,4) {};
\node (d14) at (-1,4) {};
\node (d23) at (1,4) {};
\node (d24) at (3,4) {};

\foreach \from/\to in {a1/b2,a1/b3,a1/b4, a2/b1,a2/b3,a2/b4,a3/b1,a3/b2,a3/b3,a3/b4,a4/b1,a4/b2,a4/b3,a4/b4,d13/a1,d13/a3,d14/a1,d14/a4,d23/a2,d23/a3,d24/a2,d24/a4}
\draw [-] (\from) -- (\to);
\tikzstyle{every node}=[]
\draw[left] (a1) node            {$a_1$};
\draw[right] (a2) node            {$a_2$};
\draw[left] (a3) node            {$a_3$};
\draw[right] (a4) node            {$a_4$};
\draw[below] (b1) node            {$b_1$};
\draw[below] (b2) node            {$b_2$};
\draw[below] (b3) node            {$b_3$};
\draw[below] (b4) node            {$b_4$};
\draw[above] (d13) node            {$d_{1,3}$};
\draw[above] (d14) node            {$d_{1,4}$};
\draw[above] (d23) node            {$d_{2,3}$};
\draw[above] (d24) node            {$d_{2,4}$};

\end{tikzpicture}
\caption{For the proof of \ref{442}, step (2).} \label{fig:442+}
\end{figure}
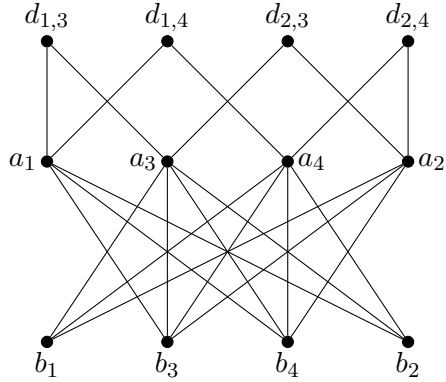
\noindent
(2) {\em There is a vertex $b_5$ not in $\{b_1\LL b_4\}$ adjacent to $a_1,a_2$}
\\
\\
Suppose not. Let $d$ be the neighbour of $a_1$ different from $d_{1,3}, d_{1,4}, b_2,b_3,b_4$. Thus $d$ is adjacent to $a_1$
and to none of $a_2,a_3,a_4$. By \ref{mate}, for $2\le i\le 4$ there exists $c_i$ different from $a_1$ that is adjacent to both
$d,b_i$. Also by \ref{mate}, there is a vertex $f$ different from $a_3$ that is adjacent to both $b_2,d_{2,3}$. So $c_2,c_3,c_4,f$
are all in $A\setminus \{a_1\LL a_4\}$. But then 
$$\{\{c_2,c_4,f,b_2\},\{b_3\},\{a_2,d_{2,3}\},\{a_3,b_1\},\{a_4,d_{2,4}\},\{b_4\}\}$$
is a 6-cluster. This proves (2).

\bigskip

By \ref{mate}, there is a vertex $c_1$
different from $a_2$ and adjacent to $b_5,b_1$; and so $c_1\ne a_1\LL a_4$. Similar there exists $c_2\ne a_1\LL a_4$
adjacent to $b_5,b_2$.
By \ref{mate}, there is a vertex $f_{1,3}$ different from $a_3$ that is adjacent to both $b_1, d_{1,3}$; and similarly there exists
$f_{1,4}$ different from $a_4$ that is adjacent to both $b_1, d_{1,4}$; there exists 
$f_{2,3}$ different from $a_3$ that is adjacent to both $b_2, d_{2,3}$; and there exists
$f_{2,4}$ different from $a_4$ that is adjacent to both $b_2, d_{1,4}$. Consequently none of $f_{1,3},f_{1,4},f_{2,3},f_{2,4}$
belongs to $\{a_1\LL a_4\}$. If they are all equal, equal to $a_5$ say, then 
$a_5$ is adjacent to $b_1,b_2,d_{1,3},d_{1,4},d_{2,3},d_{2,4}$; but then 
$$\{\{c_1,c_2,b_1,b_2, b_5\}, \{a_5\},\{a_1,d_{1,4}\},\{a_2,d_{2,3},b_3\},\{a_3,d_{1,3}\},\{a_4,d_{2,4},b_4\}\}$$
is a 6-cluster. Thus they are not all equal, and so there exist $i,j\in \{3,4\}$ such that $f_{1,i}$ is different from $f_{2,j}$.
\begin{figure}[H]
\centering

\begin{tikzpicture}[scale=0.8,auto=left]
\tikzstyle{every node}=[inner sep=1.5pt, fill=black,circle,draw]
\node (a1) at (-3,2) {};
\node (a2) at (3,2) {};
\node (a3) at (-1,2) {};
\node (a4) at (1,2) {};
\node (b1) at (-4,-1) {};
\node (b2) at (4,-1) {};
\node (b3) at (-2,-1) {};
\node (b4) at (2,-1) {};
\node (b5) at (0,-1) {};
\node (d13) at (-3,4) {};
\node (d14) at (-1,4) {};
\node (d23) at (1,4) {};
\node (d24) at (3,4) {};
\node (c1) at (-2,-2) {};
\node (c2) at (2,-2) {};
\node (f1i) at (-5,2) {};
\node (f2j) at (5,2) {};

\foreach \from/\to in {a1/b2,a1/b3,a1/b4, a2/b1,a2/b3,a2/b4,a3/b1,a3/b2,a3/b3,a3/b4,a4/b1,a4/b2,a4/b3,a4/b4,d13/a1,d13/a3,d14/a1,d14/a4,d23/a2,d23/a3,d24/a2,d24/a4, a1/b5,a2/b5, c1/b1,c1/b5,c2/b2,c2/b5, f1i/b1,f2j/b2}
\draw [-] (\from) -- (\to);
\draw (f1i) to (-2,4);
\draw (f2j) to (2,4);
\tikzstyle{every node}=[]
\draw[left] (a1) node            {$a_1$};
\draw[right] (a2) node            {$a_2$};
\draw[left] (a3) node            {$a_3$};
\draw[right] (a4) node            {$a_4$};
\draw[below] (b1) node            {$b_1$};
\draw[below] (b2) node            {$b_2$};
\draw[below] (b3) node            {$b_3$};
\draw[below] (b4) node            {$b_4$};
\draw[below] (b5) node            {$b_5$};
\draw[above] (d13) node            {$d_{1,3}$};
\draw[above] (d14) node            {$d_{1,4}$};
\draw[above] (d23) node            {$d_{2,3}$};
\draw[above] (d24) node            {$d_{2,4}$};
\draw[below] (c1) node            {$c_1$};
\draw[below] (c2) node            {$c_2$};
\draw[left] (f1i) node            {$f_{1,i}$};
\draw[right] (f2j) node            {$f_{2,j}$};

\end{tikzpicture}
\caption{For the last part of the proof of \ref{442}. $f_{1,i}$ is adjacent to $d_{1,i}$ for some $i\in \{3,4\}$, and similarly for $f_{2,j}$. The vertices $c_1,c_2$ might be equal, but $f_{1,i}$ is different from $f_{2,j}$.} \label{fig:442last}
\end{figure}
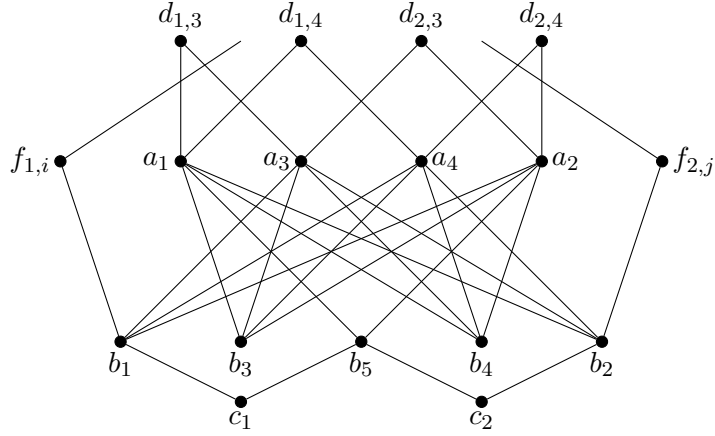
\noindent
(3) {\em There are two disjoint subsets $X,Y$ of $\{b_1,f_{1,i},b_2,f_{2,j}, c_1,c_2,b_5\}$,
both inducing connected subgraphs, with $b_1,f_{1,3}\in X$ and $b_2,f_{2,j}\in Y$, such that there is an edge between $X,Y$.}
\\
\\
If $f_{1,i}$ is adjacent to $b_2$ we may take $X=\{b_1,f_{1,i}\}$ and $Y=\{b_2,f_{2,j}\}$, so we assume that 
$f_{1,i}$ is nonadjacent to $b_2$, and similarly $f_{2,j}$ is nonadjacent to $b_1$. It follows that $f_{1,i}\ne c_2$
and $f_{2,j}\ne c_1$. So the only possible equalities between two of $f_{1,i}, f_{2,j}, c_1,c_2$ are
$f_{1,i}=c_1$, $f_{2,j}=c_2$, and $c_1=c_2$. If $f_{1,i}=c_1$, then $c_1$ is different from $f_{2,j},c_2$ and we may
set $X=\{b_1,c_1,b_5\}$ and $Y=\{b_2,f_{2,j},c_2\}$, so we assume $f_{1,i}\ne c_1$, and similarly $f_{2,j}\ne c_2$. 
But then we may set $X=\{b_1,f_{1,i}\}$ and $Y=\{b_5, c_1,c_2, b_2,f_{2,j}\}$. This proves (3).

\bigskip

But then 
$$\{X, Y,\{a_1, d_{1,3},d_{1,4}\},\{a_2,b_4,d_{2,3},d_{2,4}, b_4\},\{a_3,b_3\},\{a_4\}\}$$
is a 6-cluster. This proves \ref{442}.~\bbox

\section{Excluding $K(2,4,0)$-subgraphs}

We begin with:
\begin{thm}\label{341}
Let $G$ be a minimal candidate with bipartition $(A,B)$. Then $G$ has no $K(3,4,1)$-subgraph.
\end{thm}
\Proof
Suppose that $a_1, a_2, a_3\in A$ are all adjacent to each of $b_1\LL b_4\in B$ except $a_1b_1$.
We observe:
\\
\\
(1) {\em No vertex in $A\setminus \{a_1,a_2,a_3\}$ is adjacent to $b_1$ and to two of $b_2,b_3,b_4$.
No vertex in $B\setminus \{b_1\LL b_4\}$ is adjacent to both $a_2,a_3$. At most one vertex in $B\setminus \{b_1\LL b_4\}$ is 
adjacent to both $a_1,a_2$, and at most one is adjacent to
$a_1,a_3$,}
\\
\\
The first is because there is no $K(4,4,2)$-subgraph by \ref{442}, and the other three because there is no $K(2,5,0)$-subgraph by \ref{250}. This proves (1).
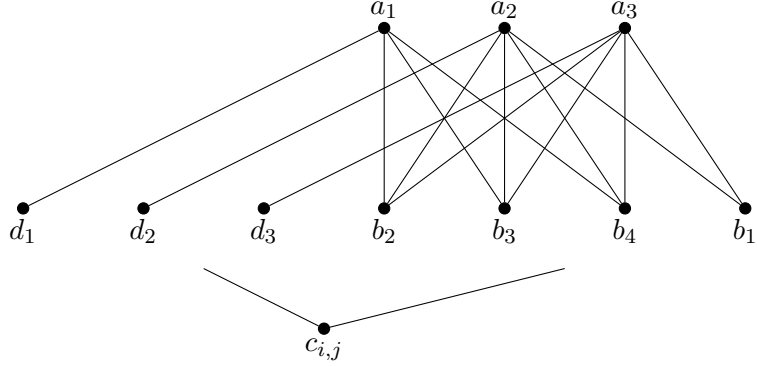
\begin{figure}[H]
\centering

\begin{tikzpicture}[scale=0.8,auto=left]
\tikzstyle{every node}=[inner sep=1.5pt, fill=black,circle,draw]
\node (a1) at (-2,3) {};
\node (a2) at (-0,3) {};
\node (a3) at (2,3) {};
\node (b1) at (4,0) {};
\node (b2) at (-2,0) {};
\node (b3) at (0,0) {};
\node (b4) at (2,0) {};
\node (d1) at (-8,0) {};
\node (d2) at (-6,0) {};
\node (d3) at (-4,0) {};
\node (c) at (-3,-2) {};

\foreach \from/\to in {a1/b2,a1/b3,a1/b4, a2/b1,a2/b2,a2/b3,a2/b4,a3/b1,a3/b2,a3/b3,a3/b4, a1/d1,a2/d2,a3/d3}
\draw [-] (\from) -- (\to);

\draw[-] (c) to (-5,-1);
\draw[-] (c) to (1,-1);
\tikzstyle{every node}=[]
\draw[above] (a1) node            {$a_1$};
\draw[above] (a2) node            {$a_2$};
\draw[above] (a3) node            {$a_3$};
\draw[below] (b1) node            {$b_1$};
\draw[below] (b2) node            {$b_2$};
\draw[below] (b3) node            {$b_3$};
\draw[below] (b4) node            {$b_4$};
\draw[below] (d1) node            {$d_1$};
\draw[below] (d2) node            {$d_2$};
\draw[below] (d3) node            {$d_3$};
\draw[below] (c) node            {$c_{i,j}$};

\end{tikzpicture}
\caption{The start of the proof of \ref{341}. $c_{i,j}$ is adjacent to $d_i$ and $b_j$.} \label{fig:341}
\end{figure}

Let $B_0$ be the set of vertices that have exactly one neighbour
in $\{a_1,a_2,a_3\}$. Thus $|B_0|$ is odd.
From (1) it follows that each of $a_1,a_2,a_3$ has a neighbour in $B_0$; let us call these neighbours
$d_1,d_2,d_3$ respectively.
For $i = 1,2,3$ and $j=1,2,3,4$ (except when $(i,j)=(1,1)$), by \ref{mate} there is a vertex $c_{i,j}$ different from $a_i$ and adjacent to $d_i, b_j$.
Hence $c_{i,j}\ne a_1, a_2, a_3$. For $i = 2,3$ choose the set $\{c_{i,1},c_{i,2},c_{i,3},c_{i,4}\}$ minimal, and choose the 
set $\{c_{1,2},c_{1,3},c_{1,4}\}$ minimal. 
\\
\\
(2) {\em Some vertex $a_4\in A\setminus \{a_1,a_2,a_3\}$ has more than one neighbour in $\{b_2,b_3,b_4\}$.}
\\
\\
Suppose not. Thus $c_{2,2},c_{2,3},c_{2,4}$ are all different.
We chose $\{c_{2,1},c_{2,2},c_{2,3},c_{2,4}\}$ minimal; so we may assume that 
either 
\begin{itemize}
\item $c_{2,1},c_{2,2},c_{2,3},c_{4,4}$ are all distinct, and each has only one neighbour in $\{b_1\LL b_4\}$; or
\item $c_{2,1}=c_{2,2}$. In this case $c_{2,1}$ is nonadjacent to $b_3,b_4$ by (1).
\end{itemize}
In both cases neither of $c_{2,1},c_{2,2}$ have a neighbour in $\{b_3,b_4\}$, and so $\{c_{1,3},c_{1,4}\}$ is disjoint from
$\{c_{2,1},c_{2,2}\}$.
But then 
$$\{\{c_{2,1},c_{2,2},d_2\},\{a_2\},\{c_{1,3}, c_{2,3},d_1,b_3\},\{c_{1,4},c_{2,4},b_4\},\{a_3,b_1\},\{a_1,b_2\}\}$$
is a 6-cluster. This proves (2).

\bigskip

We assume that $a_4$ is adjacent to $b_2,b_3$ (and possibly to $b_4$, but not to $b_1$, by (1)).
\\
\\
(3) {\em $a_4$ is adjacent to $d_2,d_3$.}
\\
\\
Suppose that $a_4,d_2$ are nonadjacent, say. Then 
$$\{\{c_{2,1},c_{2,2},c_{2,3},c_{2,4},d_2\},\{a_2\},\{b_2\},\{a_4,b_4\},\{a_3,b_1\},\{a_1,b_4\}\}$$
is a 6-cluster. This proves (3).

\begin{figure}[H]
\centering

\begin{tikzpicture}[scale=0.8,auto=left]
\tikzstyle{every node}=[inner sep=1.5pt, fill=black,circle,draw]
\node (a1) at (-2,3) {};
\node (a2) at (-0,3) {};
\node (a3) at (2,3) {};
\node (b1) at (4,0) {};
\node (b2) at (-2,0) {};
\node (b3) at (0,0) {};
\node (b4) at (2,0) {};
\node (d1) at (-8,0) {};
\node (d2) at (-6,0) {};
\node (d3) at (-4,0) {};
\node (a4) at (-3,-1.5) {};

\foreach \from/\to in {a1/b2,a1/b3,a1/b4, a2/b1,a2/b2,a2/b3,a2/b4,a3/b1,a3/b2,a3/b3,a3/b4, a1/d1,a2/d2,a3/d3, a4/b2,a4/b3,a4/d2,a4/d3}
\draw [-] (\from) -- (\to);

\tikzstyle{every node}=[]
\draw[above] (a1) node            {$a_1$};
\draw[above] (a2) node            {$a_2$};
\draw[above] (a3) node            {$a_3$};
\draw[right] (b1) node            {$b_1$};
\draw[left] (b2) node            {$b_2$};
\draw[right] (b3) node            {$b_3$};
\draw[right] (b4) node            {$b_4$};
\draw[left] (d1) node            {$d_1$};
\draw[left] (d2) node            {$d_2$};
\draw[left] (d3) node            {$d_3$};
\draw[below] (a4) node            {$a_4$};

\end{tikzpicture}
\caption{For the proof of \ref{341}, step (4).} \label{fig:341+}
\end{figure}
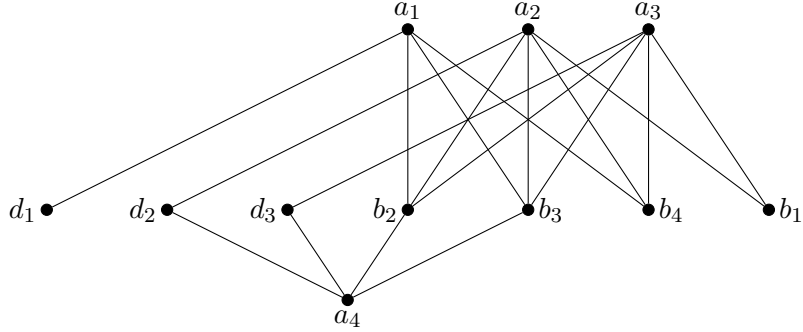
\noindent
(4) {\em If $a_4$ is nonadjacent to $b_4$, then no vertex in $A\setminus \{a_1\LL a_4\}$ has a neighbour in $\{b_2,b_3\}$
and a neighbour in $\{d_2,d_3\}$.}
\\
\\
Suppose that some $a_5 \in A\setminus \{a_1\LL a_4\}$ is adjacent to $b_2,d_2$ say. By (1), $c_{2,1}\ne a_4$, and since $a_4,b_4$
are nonadjacent, $c_{2,4}\ne a_4$. Hence
$$\{\{a_2\},\{c_{2,1},c_{2,4},a_5, d_2\},\{a_1,b_4\},\{a_3,b_1\}, \{b_2\}, \{a_4,b_3\}\}$$
is a 6-cluster. This proves (4).
\\
\\
(5) {\em $a_4$ is adjacent to $d_1$ and hence to every vertex in $B_0$.}
\\
\\
Suppose that $a_4$ is nonadjacent to $d_1$. Thus $c_{1,2}, c_{1,3}, c_{1,4}\ne a_4$. 
If $a_4$ is adjacent to $b_4$, then $c_{2,4}=a_4$ by (3) (with $b_3,b_4$ exchanged),
and at most one of $c_{1,2}, c_{1,3}, c_{1,4}=c_{2,1}$ by (1), and so we may assume that $\{c_{2,1}, c_{2,4}\}$ is disjoint from 
$\{c_{1,2}, c_{1,3}\}$.
If $a_4$ is nonadjacent to $b_4$, then by (3), neither of $c_{2,1}, c_{2,4}$ has a neighbour in $\{b_2,b_3\}$,
and so again $c_{2,1}, c_{2,4}\ne c_{1,2}, c_{1,3}$. In either case
$$\{\{b_2\},\{b_3,c_{1,2}, c_{1,3}, d_1\},\{d_2,c_{2,1}, c_{2,4}, a_4\},\{a_1,b_4\},\{a_2\},\{a_3,b_1\}\}$$
is a 6-cluster. Hence $a_4$ is adjacent to each of $d_1,d_2,d_3$, and hence to every vertex in $B_0$. This proves (5).
\\
\\
(6) {\em $a_4$ is adjacent to $b_4$.}
\\
\\
Suppose not. Since $a_4$ has degree six and is nonadjacent to $b_1$ by (1), it has a neighbour $b_5\ne b_1\LL b_4,d_1,d_2,d_3$. By \ref{mate},
for each $v\in \{b_2,b_3,d_1,d_2,d_3\}$ there is a vertex $c(v)$ different from $a_4$ adjacent to $b_5, v$; and hence
$c(v)\ne a_1\LL a_4$. But then 
$$\{\{a_4\},\{b_5,c(b_2), c(b_3), c(d_1), c(d_2), c(d_3)\},\{a_1,d_1,b_4\},\{a_2,b_1,d_2\},\{a_3,d_3\},\{b_2\}\}$$
is a 6-cluster. This proves (6).

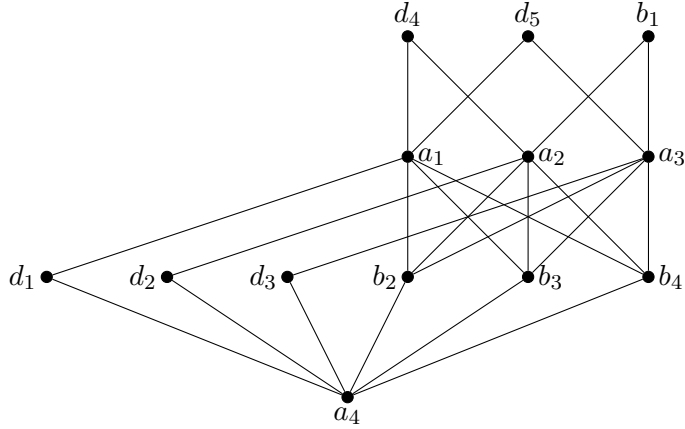
\begin{figure}[H]
\centering

\begin{tikzpicture}[scale=0.8,auto=left]
\tikzstyle{every node}=[inner sep=1.5pt, fill=black,circle,draw]
\node (a1) at (-2,3) {};
\node (a2) at (-0,3) {};
\node (a3) at (2,3) {};
\node (b1) at (2,5) {};
\node (b2) at (-2,1) {};
\node (b3) at (0,1) {};
\node (b4) at (2,1) {};
\node (d1) at (-8,1) {};
\node (d2) at (-6,1) {};
\node (d3) at (-4,1) {};
\node (a4) at (-3,-1) {};
\node (d4) at (-2,5) {};
\node (d5) at (0,5) {};

\foreach \from/\to in {a1/b2,a1/b3,a1/b4, a2/b1,a2/b2,a2/b3,a2/b4,a3/b1,a3/b2,a3/b3,a3/b4, a1/d1,a2/d2,a3/d3, a4/b2,a4/b3,a4/d2,a4/d3, a4/d1,a4/b4, d4/a1,d4/a2,d5/a1,d5/a3}
\draw [-] (\from) -- (\to);

\tikzstyle{every node}=[]
\draw[right] (a1) node            {$a_1$};
\draw[right] (a2) node            {$a_2$};
\draw[right] (a3) node            {$a_3$};
\draw[above] (b1) node            {$b_1$};
\draw[left] (b2) node            {$b_2$};
\draw[right] (b3) node            {$b_3$};
\draw[right] (b4) node            {$b_4$};
\draw[left] (d1) node            {$d_1$};
\draw[left] (d2) node            {$d_2$};
\draw[left] (d3) node            {$d_3$};
\draw[below] (a4) node            {$a_4$};
\draw[above] (d4) node            {$d_4$};
\draw[above] (d5) node            {$d_5$};

\end{tikzpicture}
\caption{For the last part of the proof of \ref{341}.} \label{fig:341last}
\end{figure}

\bigskip

Since $a_4$ has degree only six, it follows from (5) that $|B_0|\le 4$. Since $|B_0|$ is odd, it follows that $B_0=\{d_1,d_2,d_3\}$, 
and there are two vertices $d_4,d_5\in B\setminus \{b_1\LL b_4\}$ that have two neighbours in $\{a_1,a_2,a_3\}$. By (1) we may assume that $d_4, d_5$ 
are adjacent to $a_1,a_2$ and to $a_1,a_3$ respectively.

We do not obtain a smaller candidate by contracting $\{d_4,a_1,b_2\}, \{d_5,a_3,b_3\},\{b_1,a_2,b_4\}$ into $B$;
and so there is a vertex $a_5\ne a_1\LL a_4$, adjacent either to both $b_2,d_4$, or to both $b_3,b_5$, or to both $b_1,b_4$.
There is symmetry between $b_1,d_4,d_5$ (see figure \ref{fig:341last}), so we may assume that $a_5$ is adjacent to $b_1,b_4$. By \ref{mate},
there is a vertex $c(d_2)$ different from $a_2$ adjacent to both $b_1,d_2$; 
a vertex $c(d_3)$ different from $a_3$ adjacent to both $b_1,d_3$; and a vertex $c(d_5)$ different from $a_3$
adjacent to both $b_1,d_5$. It follows that $c(d_2), c(d_3)\ne a_1\LL a_4$. But then 
$$\{\{b_1,a_5,c(d_2),c(d_3), c(d_5), d_2\},\{a_2,b_3\},\{a_1,d_5\},\{a_3,d_3\},\{a_4,b_2\},\{b_4\}\}$$
is a 6-cluster. This proves \ref{341}.~\bbox

\begin{thm}\label{240}
Let $G$ be a minimal candidate with bipartition $(A,B)$. Then $G$ has no $K(2,4,0)$-subgraph.
\end{thm}
\Proof
Suppose that $a_1, a_2\in A$ are both adjacent to each of $b_1\LL b_4\in B$. No other vertex is adjacent to both $a_1,a_2$, since there
is no $K(2,5,0)$-subgraph by \ref{250}. Let $b_5,b_6,b_7,b_8$ be the vertices in $B$ that have exactly one neighbour in 
$\{a_1,a_2\}$, where $b_5,b_6$ are adjacent to $a_1$, and $b_7,b_8$ to $a_2$.

No other vertex is adjacent to more than two of of $b_1\LL b_4$, since there is no 
$K(3,4,0)$- or $K(3,4,1)$-subgraph by \ref{340} and \ref{341}. Let $H$ be the cover graph with respect to $a_1,a_2,b_1,b_2,b_3,b_4$;
then $\chi(H)\ge 3$ by \ref{colouring}, and so $H$ has a triangle, say with vertices $b_1,b_2,b_3$. Consequently there are three
vertices $c_1,c_2,c_3$, such that $c_i$ is adjacent to the two vertices in $\{b_1,b_2,b_3\}\setminus \{b_i\}$, and $c_i$
is nonadjacent to $b_i$.  
\begin{figure}[H]
\centering

\begin{tikzpicture}[scale=0.8,auto=left]
\tikzstyle{every node}=[inner sep=1.5pt, fill=black,circle,draw]
\node (a1) at (-2.5,2) {};
\node (a2) at (.5,2) {};
\node (b1) at (-3,0) {};
\node (b2) at (-1,0) {};
\node (b3) at (1,0) {};
\node (b4) at (-1,4) {};
\node (b5) at (-7,0) {};
\node (b6) at (-5,0) {};
\node (b7) at (3,0) {};
\node (b8) at (5,0) {};
\node (c1) at (1,-2) {};
\node (c2) at (-1,-2) {};
\node (c3) at (-3,-2) {};

\foreach \from/\to in {a1/b1,a1/b2,a1/b3,a1/b4, a2/b1,a2/b2,a2/b3,a2/b4, a1/b5,a1/b6,a2/b7,a2/b8, c1/b2,c1/b3,c2/b1,c2/b3,c3/b1,c3/b2}
\draw [-] (\from) -- (\to);
\tikzstyle{every node}=[]
\draw[above left] (a1) node            {$a_1$};
\draw[above right] (a2) node            {$a_2$};

\draw[below left] (b1) node            {$b_1$};
\draw[below] (b2) node            {$b_2$};
\draw[below right] (b3) node            {$b_3$};
\draw[above] (b4) node            {$b_4$};
\draw[below] (b5) node            {$b_5$};
\draw[below] (b6) node            {$b_6$};
\draw[below] (b7) node            {$b_7$};
\draw[below] (b8) node            {$b_8$};
\draw[below] (c1) node            {$c_1$};
\draw[below] (c2) node            {$c_2$};
\draw[below] (c3) node            {$c_3$};

\end{tikzpicture}
\caption{For the proof of \ref{240}.} \label{fig:240}
\end{figure}
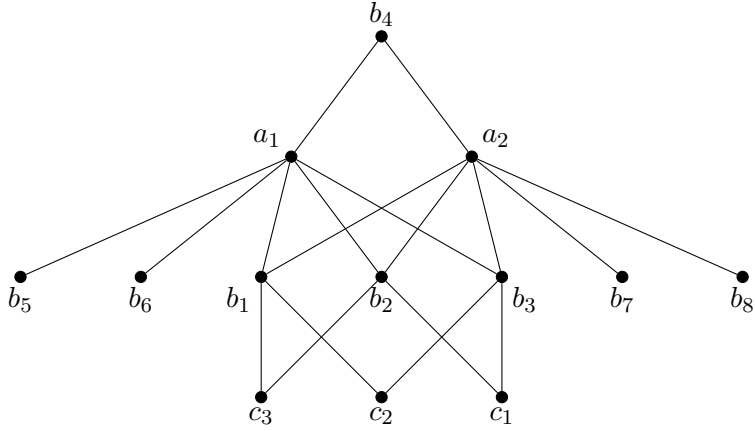

Let $Z=\{a_1,a_2,b_1,b_2,b_3,b_4, c_1,c_2,c_3\}$, and let
$\mathcal{C}$ be the set of components of $G\setminus Z$ that contain at least one of $b_5,b_6,b_7,b_8$. For each $C\in \mathcal{C}$,
let $N(C)$ be the set of vertices in $\{b_1,b_2,b_3, c_1,c_2,c_3\}$ that have a neighbour in $V(C)$. 
\\
\\
(1) {\em If $C\in \mathcal{C}$, then $b_4$ has a neighbour in $V(C)$, and either
\begin{itemize}
\item $N(C)= \{b_i, c_i\}$ for some $i\in \{1,2,3\}$; or
\item $N(C)= \{c_i,c_j\}$ for some two distinct $i,j\in \{1,2,3\}$.
\end{itemize}
}
\noindent 
Let $b_5\in V(C)$ say. By \ref{mate}, for $i = 1,2,3,4$ there is a vertex  different
from $a_1$ adjacent to both $b_5, b_i$. Since none of $c_1,c_2,c_3$ is adjacent to $b_4$, it follows that $b_4$ has a neighbour in $V(C)$.
Also, for $1\le i\le 3$, either $b_i\in N(C)$ or $N(C)$ contains a vertex in
$\{c_1,c_2,c_3\}\setminus \{c_i\}$. In summary, $N(C)$ contains a member of each of the sets 
$\{b_1,c_2,c_3\}$, $\{b_2,c_3,c_1\}$, $\{b_3,c_1,c_2\}$. 
Consequently, if $|N(C)|=2$ then 
the claim holds, so we assume that $|N(C)|\ge 3$. Since
$$\{\{c_3,b_1\},\{c_2,b_3\},\{c_1,b_2\},\{a_1\},\{a_2,b_4\}, V(C)\}$$
is not a 6-cluster, it follows that $N(C)$ is disjoint from one of the sets $\{c_3,b_1\},\{c_2,b_3\},\{c_1,b_2\}$, 
and we may assume from the symmetry that 
$b_1,c_3\notin N(C)$. Since $N(C)$ contains a member of
$\{b_1,c_2,c_3\}$, it follows that $c_2\in N(C)$. By a similar argument, $N(C)$ is disjoint from one of
$\{c_3,b_2\},\{c_2,b_1\},\{c_1,b_3\}$, and hence from one of $\{c_3,b_2\},\{c_1,b_3\}$. Since 
$|N(C)|\ge 3$, $N(C)\cap \{c_1,b_3\}\ne \emptyset$; and hence $b_2\notin N(C)$. Since
$|N(C)|\ge 3$ it follows that $c_1,b_3\in N(C)$. But then
$$\{\{c_2,c_3,b_1\},\{c_1,b_2\},\{b_3\},\{a_1\},\{a_2,b_4\}, V(C)\}$$
is a 6-cluster.  This proves (1).
\\
\\
(2) {\em If $C\in \mathcal{C}$, then every vertex of $V(C)\cap \{b_5,b_6,b_7,b_8\}$ is adjacent to every vertex of $N(C)\cap \{c_1,c_2,c_3\}$.}
\\
\\
We may assume that $b_5\in V(C)$ and $c_3\in N(C)$, and we must show that $b_5,c_3$
are adjacent. It follows from (1) that $b_1,b_2\notin N(C)$. But
for $i = 1,2$ there is a vertex  different
from $a_1$ adjacent to both $b_5, b_i$, and this vertex is one of $c_1,c_2,c_3$ since $b_1,b_2\notin N(C)$. Thus $b_5$ is adjacent to
one of $c_1,c_3$, and also to one of $c_2,c_3$. But $b_5$ is not adjacent to both $c_1,c_2$ since not all $c_1,c_2,c_3\in N(C)$; 
and so $b_5$ is adjacent to $c_3$. This proves (2).
\\
\\
(3) {\em Let $C\in \mathcal{C}$ contain one of $b_5,b_6,b_7,b_8$; then $N(C)\cap \{b_1,b_2,b_3\}=\emptyset$.}
\\
\\
Suppose that $b_1, b_5\in V(C)$ say. Thus $N(C)$ contains $c_1$ and none of $b_2,b_3,c_2,c_3$. 
Let $C'\in \mathcal{C}$ contain $b_6$.
We claim that $b_6$ is also
adjacent to $c_1$; because suppose not. Then $b_6\notin V(C)$, by (2), and so $C'\ne C$. Since
$b_6$ is nonadjacent to $c_1$, (2) implies that $c_1\notin N(C')$. 
Since $b_5,b_6$
are both adjacent to $a_1$, they have a second common neighbour $c$ say. But $c\notin \{c_1,c_2,c_3\}$, since $c_1$ is not adjacent to $b_6$, 
and $c_2,c_3\notin N(C)$ and so are not adjacent to $b_5$. Also $c\ne a_1,a_2$; so $c\notin Z$, contradicting that $C\ne C'$
are components of $G\setminus Z$. This proves that $b_6$ is adjacent to $c_1$. 

We claim also that $c_2,c_3\notin N(C')$. There is a symmetry exchanging $c_2,c_3$ and fixing $c_1$, so we suppose without loss of generality
that $c_2\in N(C')$, and hence $C\ne C'$. But then
$$\{\{V(C')\cup \{c_2\},V(C)\cup \{c_1,b_3\}, \{b_1\},\{b_2,c_3\}, \{a_1\},\{a_2,b_4\}\}$$
is a 6-cluster. This proves that  $c_2,c_3\notin N(C')$, and hence $N(C')=N(C)$.

Now $a_1,c_1$ have four common neighbours,
namely $b_2,b_3,b_5,b_6$; and so they satisfy the same conditions as $a_1,a_2$. At the start of this proof, we showed the existence
of $c_1,c_2,c_3$, that with three of $b_1,b_2,b_3$ induce a 6-cycle. Consequently the same is true
for $b_2,b_3,b_5,b_6$, and in particular there are two vertices $c_4,c_5$, different from $a_1,c_1$, 
that have a neighbour in $\{b_2,b_3\}$. Since $c_2,c_3\notin N(C)= N(C')$, it follows that $c_4,c_5\notin Z$, and hence
at least one of $b_2,b_3$ belongs to $N(C)= N(C')$, contrary to (1). This proves (3).

\bigskip

For $i = 5,6,7,8$, let $C_i\in \mathcal{C}$ contain $b_i$ (they are not necessarily all different).
From (1), (2) and (3) it follows that each of $b_5,b_6,b_7,b_8$ is adjacent to two of $c_1,c_2,c_3$. They are not all
adjacent to the same two of $c_1,c_2,c_3$, since there is no $K(2,5,0)$-subgraph by \ref{250}; so there exists 
$C\in \{C_5,C_6\}$ and $C'\in \{C_7,C_8\}$ with $N(C)\ne N(C')$. Thus we may assume that
$c_2,c_3\in N(C_5)$ and $c_1,c_2\in N(C_7)$. In particular $C_5\ne C_7$; and by (2), $b_5, b_7$ are adjacent to $c_2$. By \ref{mate}
they have another common neighbour, say $c$. Thus $c\notin \{a_1,a_2\}$ because $b_5, b_7$ each have only one and different neighbours
in that set; $c\notin \{c_1,c_2,c_3\}$, since $c\ne c_2$ from its definition, and $c_1\notin N(C_5)$, and $c_3\notin N(C_7)$;
and $c\notin V(G)\setminus Z$ since $C\ne C'$, a contradiction. This proves \ref{240}.~\bbox

\section{The end}

Next we need the following lemma. Let $H$ be a complete graph with six vertices. 
If $C$ is a triangle of $H$, a {\em $C$-path} means a path of $H$ with vertex set $C$. Let $C_1\LL C_k$ be triangles
of $H$, not            
necessarily all different. We denote by $M(C_1\LL C_k)$ the graph with vertex set $V(H)$ and edge set the set of all edges $uv$ 
of $H$ such that $\{u,v\}$ is not a subset of any of $C_1\LL C_k$. Let $J$ be the graph obtained from a six-vertex complete graph
by deleting four edges, the edges of two disjoint three-vertex paths. We observe that $J$ admits
a 5-cluster (one of the five sets contains the two vertices of degree three, the others are singletons). 

\begin{thm}\label{game}
Let $H, C_1\LL C_k$ be as above.
Then for $1\le i\le k$ there is a $C_i$-path $P_i$, such that 
$$M(C_1\LL C_k)\cup P_1\cupcup P_k$$ 
has a subgraph isomorphic to $J$.
\end{thm}
\Proof
We proceed by induction on $k$. Suppose first that $u,v\in V(H)$ and two of $C_1\LL C_k$ contain $u,v$, say $C_1,C_2$.
Let $C_1=\{u,v,w\}$ say. Let $P_1$ be the path $u\DD w\DD v$. From the inductive hypothesis, for $2\le i\le k$ there is a
$C_i$-path $P_i$ such that $M(C_2\LL C_k)\cup P_2\cupcup P_k$ admits a 5-cluster. But every edge of $M(C_2\LL C_k)$
that is not an edge of $M(C_1\LL C_k)$ is one of $uw,vw$, and they are edges of $P_1$; so 
$M(C_1\LL C_k)\cup P_1\cupcup P_k$ contains a copy of $J$ as required.

So we may assume that no two of $C_1\LL C_k$ share more than one vertex. 
Every subgraph of $H$ obtained by deleting at most two 
edges contains a copy of $J$, so we may assume that $k\ge 3$ and hence no two of $C_1\LL C_k$ are disjoint 
(because if $C_1\cap C_2=\emptyset$ then $C_3$ shares two vertices with one of them).
Let $V(H)=\{h_1\LL h_6\}$; then we may assume that $C_1=\{h_1,h_2,h_3\}$, $C_2=\{h_1,h_4,h_5\}$, $C_3=\{h_2,h_4,h_6\}$ and either $k=3$, or $k=4$
and $C_4=\{h_3,h_5,h_6\}$. Define $P_1=h_2\DD h_1\DD h_3$, $P_2= h_1\DD h_5\DD h_4$, $P_3=h_4\DD h_2\DD h_6$, and if $k=4$, define $P_4=h_3\DD h_6\DD h_5$. 
If $k=4$, 
$M(C_1\LL C_k)\cup P_1\cupcup P_k$ is isomorphic to $J$, and if $k=3$, $M(C_1\LL C_k)\cup P_1\cupcup P_k$ contains a copy of $J$.
This proves \ref{game}.~\bbox

We deduce \ref{mainthm}, which we restate as follows:
\begin{thm}\label{nocandidate}
No graph is a candidate.
\end{thm}
\Proof
Assume some graph is a candidate, and let $G$ be a minimal candidate, with bipartition $(A,B)$. Let $a\in A$, and let its
neighbours be $b_1\LL b_6$. By \ref{colouring} (or by \ref{mate}), the cover graph $H$ with respect to $a,b_1\LL b_6$ is 
complete. No vertex different from $a$ has more than three neighbours in $\{b_1\LL b_6\}$, since there is no $K(2,6,0)$-, 
$K(2,5,0)$- or $K(2,4,0)$-subgraph,  by \ref{260}, \ref{250} and 
\ref{240}. Consequently there is a set of vertices $a_1\LL a_\ell\in A\setminus \{a\}$, each with two or three neighbours in 
$\{b_1\LL b_6\}$, such that for all distinct $u,v\in \{b_1\LL b_6\}$, there exists $i\in \{1\LL \ell\}$ such that $a_i$
is adjacent to $u,v$. We may assume that $a_1\LL a_k$ have three neighbours in $\{b_1\LL b_6\}$ and $a_{k+1}\LL a_\ell$ have two.
For $1\le i\le k$ let $C_i$ be the set of three neighbours of $a_i$ in $\{b_1\LL b_6\}$. By \ref{game}, for $1\le i\le k$
there is a $C_i$-path $P_i$ such that $M(C_1\LL C_k)\cup P_1\cupcup P_k$ contains a copy of $J$ as a subgraph, and hence admits a 
5-cluster $\{X_1\LL X_5\}$ say. Let $p_i$ be the middle vertex of $P_i$ for $1\le i\le k$, and let $p_i$ be one of the two 
members of $C_i$ for $k+1\le i\le \ell$.

For each $v\in \{b_1\LL b_6\}$ let $C(v)$ be the union of $\{v\}$ and all the vertices $c_i$ with $1\le i\le \ell$
such that $p_i=v$. 
For $1\le j\le 5$ let $Y_j$ be the union of the sets $C(v)$ over all $v\in X_j$.
We claim that 
$\{\{a\},Y_1\LL Y_5\}$ is a 6-cluster. To see this, we must check that these six sets are pairwise disjoint subsets of $V(G)$,
which is clear; that $Y_1\LL Y_5$ each contain a neighbour of $a$, which is true since $X_1\LL X_5$ are nonempty subsets of
$\{b_1\LL b_6\}$; that each of the sets $Y_i$ induces a connected subgraph of $G$; and that for $1\le i<i'\le 5$, some vertex in 
$Y_i$ has a neighbour in $Y_{i'}$. To see both these final statements, it suffices to show that if $u,v\in \{b_1\LL b_6\}$
are adjacent in $M(C_1\LL C_k)\cup P_1\cupcup P_k$, there is an edge of $G$ between $C(u), C(v)$. To see this, there
are two cases: $uv\in E(M(C_1\LL C_k))$, and $uv\in E(P_i)$ for some $i\in \{1\LL k\}$. In the first case, there exists $i$ with $k+1\le i\le \ell$
such that $c_i$ is adjacent to both $u,v$, and since $p_i$ is one of $u,v$, the claim holds. In the second case,
let $uv\in E(P_i)$ for some $i\in \{1\LL k\}$; then one of $u,v$ equals $p_i$, say $v$; the other, $u$, is an end of $P_i$;
and in $G$ there is an edge between $c_i$ and $u$. This proves \ref{nocandidate}.~\bbox

\end{document}